\documentclass[english]{article}
\usepackage[T1]{fontenc}
\usepackage[latin9]{inputenc}
\usepackage{geometry}
\usepackage{amsmath}
\usepackage{amssymb}
\usepackage{graphicx}
\usepackage{babel}
\begin{document}

\title{A Spectral Theory for Tensors }

\author{Edinah K. Gnang%
\thanks{Department of Computer Science, Rutgers University, Piscataway, NJ
08854-8019 USA%
}, Ahmed Elgammal%
\thanks{Department of Computer Science, Rutgers University, Piscataway, NJ
08854-8019 USA%
}, Vladimir Retakh%
\thanks{Department of Mathematics, Rutgers University, Piscataway, NJ 08854-8019
USA%
}}
\maketitle
\begin{abstract}
In this paper we propose a general spectral theory for tensors. Our
proposed factorization decomposes a tensor into a product of orthogonal
and scaling tensors. At the same time, our factorization yields an
expansion of a tensor as a summation of outer products of lower order
tensors. Our proposed factorization shows the relationship between
the eigen-objects and the generalised characteristic polynomials.
Our framework is based on a consistent multilinear algebra which explains
how to generalise the notion of matrix hermicity, matrix transpose,
and most importantly the notion of orthogonality. Our proposed factorization
for a tensor in terms of lower order tensors can be recursively applied
so as to naturally induces a spectral hierarchy for tensors. 
\end{abstract}

\section{Introduction}

In 1762 Joseph Louis Lagrange formulated what is now known as the
eigenvalue - eigenvector problem, which turns out to be of significant
importance in the understanding several phenomena in applied mathematics
as well as in optimization theory. The spectral theory for matrices
is widely used in many scientific and engineering domains.

In many scientific domains, data are presented in the form of tuples
or groups, which naturally give rise to tensors. Therefore, the generalization
of the eigenvalue-eigenvector problem for tensors is a fundamental
question with broad potential applications. Many researchers suggested
different forms of tensor decompositions to generalize the concepts
of eigenvalue-eigenvector and Singular Value Decomposition.

In this paper we propose a mathematical framework for high-order tensors
algebra based on a high-order product operator. This algebra allows
us to generalize familiar notions and operations from linear algebra
including dot product, matrix adjoints, hermicity, permutation matrices,
and most importantly the notion of orthogonality. Our principal result
is to establish a rigorous formulation of tensor spectral decomposition
through the general spectral theorem. We prove the spectral theorem
for hermitian finite order tensors with norm different from $1$.
Finally we point out that one of the fundamental consequence of the
spectral theorem is the existence of a spectral hierarchy which determines
a given hermitian tensor of finite order.

There are certain properties that a general spectral theory is expected
to satisfy. The most fundamental property one should expect from a
general formulation of the spectral theorem for tensors is a factorization
of a cubic tensor into a certain number of cubic tensors of the same
dimensions. Our proposed factorization decomposes a Hermitian tensor
into a product of orthogonal and scaling tensors. Our proposed factorization
also extends to handle non-Hermitian tensors. Furthermore our proposed
factorization offers an expansion of a tensor as a summation of lower
order tensors that are obtained through outer products. Our proposed
factorization makes an explicit connection between the eigen-objects
and the reduced set of characteristic polynomials. The proposed framework
describes the spectral hierarchy associated with a tensor. Finally
the framework aims to extend linear algebraic problems found in many
domains to higher degree algebraic formulations of corresponding problems.

The organization of this paper is as follows; Section {[}2{]} reviews
the state of the art in tensor decomposition and its relation to the
proposed formulation. Section {[}3{]} introduces our proposed tensor
algebra for order three tensors. Section {[}4{]} introduces and proves
our proposed spectral theorem for order three tensors. Section {[}5{]}
discusses some important properties following from the proposed spectral
decomposition. Section {[}6{]} proposes a computational framework
for describing the characteristic polynomials of a tensor. Section
{[}7{]} generalizes the introduced concepts to higher order tensors
and introduces the notion of the spectral hierarchy. Section {[}8{]}
discusses in details the relation between the proposed framework and
some existing tensor decomposition frameworks. Section {[}9{]} concludes
the paper with a discussion on the open directions.

\section{State of the art in tensor decomposition}

\subsection{Generalizing Concepts from Linear Algebra}

In this section we recall the commonly used notation by the multilinear
algebra community where a $k$-tensor denotes a multi-way array with
$k$ indices \cite{workshop04}. Therefore, a vector is a $1$-tensor
and a matrix is a $2$-tensor. A $3$-tensor $\boldsymbol{A}$ of
dimensions $m\times n\times p$ denotes a rectangular cuboid array
of numbers. The array consists of $m$ rows, $n$ columns, and $p$
depths with the entry $a_{i,j,k}$ occupying the position where the
$i^{th}$ row, the $j^{th}$ column, and the $k^{th}$depth meet.
For many purposes it will suffice to write

\begin{equation}
\boldsymbol{A}=\left(a_{i,j,k}\right)\;\left(1\le i\le m;\,1\le j\le n;\,1\le k\le p\right),
\end{equation}
we now introduce generalizations of complex conjugate and inner product
operators. \\
The \emph{order $p$ conjugates} of a scalar complex number $z$ are
defined by: 
\begin{equation}
z^{\mathfrak{c}_{p}^{j}}\equiv\sqrt{\Re^{2}\left(z\right)+\Im^{2}\left(z\right)}\;\exp\left\{ i\times\arctan\left\{ \frac{\Im\left(z\right)}{\Re\left(z\right)}\right\} \times\exp\left\{ i\frac{2\pi\, j}{p}\right\} \right\} 
\end{equation}
where $\Im\left(z\right)$ and $\Re\left(z\right)$ respectively refer
to the imaginary and real part of the complex number $z$, equivalently
rewritten as 
\begin{equation}
z^{\mathfrak{c}_{p}^{j}}\equiv\left|z\right|\exp\left\{ i\times\angle_{z}\times\exp\left\{ i\frac{2\pi\, j}{p}\right\} \right\} ,
\end{equation}
from which it follows that 
\begin{equation}
\left|z\right|^{p}=\prod_{1\le j\le p}z^{\mathfrak{c}_{p}^{j}}\:.
\end{equation}
The particular \emph{inner product} operator that we introduce relates
the inner product of a $p$-tuple of vectors in $\mathbb{C}^{l}$
to a particular $\ell_{p}$ norm operator $\mathbb{C}^{l}$ in a way
quite similar to the way the inner product of pairs of vectors relate
to the usual $\ell_{2}$ vector norm. We refer to the norm operator
$\left\Vert \:\right\Vert _{\ell_{p}}:\mathbb{C}^{l}\rightarrow\mathbb{R}^{+}$
(for every integer $p\ge2$) as the $\ell_{p}$ norm defined for an
arbitrary vector $\boldsymbol{x}\equiv\left(\begin{array}{ccc}
x(1), & \cdots, & x(l)\end{array}\right)\in\mathbb{C}^{l}$ by
\begin{equation}
\left\Vert \boldsymbol{x}\right\Vert _{\ell_{p}}\equiv\left[\sum_{1\le k\le l}\;\prod_{1\le j\le p}\left(x(k)\right)^{\mathfrak{c}_{p}^{p-j}}\right]^{\frac{1}{p}},
\end{equation}
the inner product operator for a $p$-tuple of vectors in $\mathbb{C}^{l}$
denoted $\left\langle \;\right\rangle :\;\left(\mathbb{C}^{l}\right)^{p}\rightarrow\mathbb{C}$
is defined by

\begin{equation}
\left\langle \boldsymbol{v}_{k}\right\rangle _{0\le k\le p}\equiv\sum_{1\le j\le l}\left\{ \prod_{0\le k\le p-1}\left(v_{k}(j)\right)^{\mathfrak{c}_{p}^{p-j}}\right\} 
\end{equation}
some of the usual properties of inner products follow from the definition 

\begin{equation}
\left\langle \left(\boldsymbol{x}_{1}+\boldsymbol{y}_{1}\right)\,;\,\boldsymbol{z}_{2}\,;\,\cdots\,;\,\boldsymbol{z}_{l}\right\rangle =\left\langle \boldsymbol{x}_{1}\,;\,\boldsymbol{z}_{2}\,;\,\cdots\,;\,\boldsymbol{z}_{l}\right\rangle +\left\langle \boldsymbol{y}_{1}\,;\,\boldsymbol{z}_{2}\,;\,\cdots\,;\,\boldsymbol{z}_{l}\right\rangle 
\end{equation}
and most importantly the fact that
\begin{equation}
\left\langle \underbrace{\boldsymbol{z}\,;\,\boldsymbol{z}\,;\,\cdots\,;\,\boldsymbol{z}\,;\,\boldsymbol{z}}_{p\: operands}\right\rangle \ge0
\end{equation}
and 
\begin{equation}
\left\langle \underbrace{\boldsymbol{z}\,;\,\boldsymbol{z}\,;\,\cdots\,;\,\boldsymbol{z}\,;\,\boldsymbol{z}}_{p\: operands}\right\rangle =0\Leftrightarrow\boldsymbol{z}=\boldsymbol{0}.
\end{equation}
We point out that the definitions of inner products is extended naturally
to tensors as illustrated bellow

\begin{equation}
\left\langle \boldsymbol{A},\boldsymbol{B}\right\rangle \equiv\sum_{1\le m,n\le l}a_{m,n}\times\left(b_{n,m}\right)^{\mathfrak{c}_{2}^{1}}
\end{equation}
\begin{equation}
\left\langle \boldsymbol{A},\boldsymbol{B},\boldsymbol{C}\right\rangle \equiv\sum_{1\le m,n,p\le l}a_{m,n,p}\times\left(b_{p,m,n}\right)^{\mathfrak{c}_{3}^{2}}\times\left(c_{n,p,m}\right)^{\mathfrak{c}_{3}^{1}},
\end{equation}
More generally for arbitrarly finite order tensor the inner product
for the family of tensors $\left\{ \boldsymbol{A}^{(t)}=\left(a_{i_{1},i_{2},\cdots,i_{n}}^{(t)}\right)\right\} _{1\le t\le n}$
is defined by: 

\begin{equation}
\left\langle \boldsymbol{A}^{(t)}\right\rangle _{1\le t\le n}\equiv\sum_{1\le i_{1},i_{2},\cdots,i_{n}\le l}\left(\prod_{0\le t\le n-1}\left(a_{i_{1+\left(t-1\right)},\cdots,i_{n+\left(t-1\right)}}^{(t)}\right)^{\mathfrak{c}_{n}^{p-t}}\right)
\end{equation}
note that the addition in the indices are performed modulo $n$.\\
Generalization of other concepts arising from linear algebra have
been investigated quite extensively in the literature. Cayley in \cite{Cayley_A.}
instigated investigations on hyperdeterminants as a generalization
of determinants. Gelfand, Kapranov and Zelevinsky followed up on Cayley's
work on the subject of hyperdeterminants by relating hyperdeterminants
to $X$-discriminants in their book \cite{Gelfand_et_al.}.\\
A recent approach for generalizing the concept of eigenvalue and eigenvector
has been proposed by Liqun Qi in \cite{Qi_L2,Qi_L1} and followed
up on by Lek-Heng Lim\cite{L.-H._Lim}, Cartwright and Sturmfels \cite{DustC_SturmB}.
The starting point for their approach will be briefly summarized using
the notation introduced in the book \cite{Gelfand_et_al.}. Assuming
a choice of a coordinate system $\boldsymbol{x}_{j}=\left(x_{j}(0),x_{j}(1),\cdots,x_{j}(k_{j})\right)$
associated with each one of the vector space $V_{j}\equiv\left(\mathbb{R}^{+}\right)^{k_{j}+1}$.
We consider a multilinear function $f\::\;\bigotimes_{t=1}^{r}V_{t}\rightarrow\mathbb{R}^{+}$
expressed by :

\begin{equation}
f\left(\boldsymbol{x}_{1},\boldsymbol{x}_{2},\cdots,\boldsymbol{x}_{r-1},\boldsymbol{x}_{r}\right)=\sum_{i_{1},\cdots,i_{r}}a_{i_{1},\cdots,i_{r}}x_{1}(i_{1})\cdots x_{r}(i_{r}),
\end{equation}
equivalently the expression above can be rewritten as

\begin{equation}
f\left(\boldsymbol{x}_{1},\boldsymbol{x}_{2},\cdots,\boldsymbol{x}_{r-1},\boldsymbol{x}_{r}\right)\equiv\left\langle \boldsymbol{x}_{1},\boldsymbol{x}_{2},\cdots,\boldsymbol{x}_{r-1},\boldsymbol{x}_{r}\right\rangle _{\boldsymbol{A}}.
\end{equation}
which of course is a natural generalization of bilinear forms associated
with a matrix representation of a linear map for some choice of coordinate
system 

\begin{equation}
f\left(\boldsymbol{x}_{1},\boldsymbol{x}_{2}\right)=\sum_{i_{1},i_{2}}a_{i_{1},i_{2}}\: x_{1}(i_{1})\, x_{2}(i_{2})\equiv\left(\boldsymbol{x}_{1}\right)^{T}\boldsymbol{A}\:\boldsymbol{x}_{1}\equiv\left\langle \boldsymbol{x}_{1},\boldsymbol{x}_{2}\right\rangle _{\boldsymbol{A}}.
\end{equation}
It follows from the definition of the multilinear function $f$ that
the function induces $r$ not necessarily distinct multilinear projective
maps denoted by $f_{k}:\;\bigotimes_{\begin{array}{c}
t=1\\
t\ne k
\end{array}}^{r}V_{t}\rightarrow V_{k}$ expressed as :

\begin{equation}
f_{k}\left(\boldsymbol{x}_{1},\boldsymbol{x}_{2},\cdots\boldsymbol{x}_{k-1},\boldsymbol{x}_{k+1},\cdots,\boldsymbol{x}_{r}\right)=\sum_{i_{1},\cdots,i_{k-1},i_{k+1}\cdots,i_{r}}a_{i_{1},\cdots,i_{r}}x_{1}(i_{1})\, x_{2}(i_{2})\cdots\: x_{k-1}(i_{k-1})\: x_{k+1}(i_{k+1})\cdots\, x_{r}(i_{r})
\end{equation}
The various formulations of eigenvalue eigenvector problems as proposed
and studied in \cite{Qi_L2,Qi_L1,DustC_SturmB,L.-H._Lim} arise from
investigating solutions to equations of the form:

\begin{equation}
f_{k}\left(\boldsymbol{x},\cdots,\boldsymbol{x}\right)=\lambda\cdot\boldsymbol{x}
\end{equation}

Applying symmetry arguments to the tensor $\boldsymbol{A}$ greatly
reduces the number of map $f_{k}$ induced by $\boldsymbol{A}$. For
instance if $\boldsymbol{A}$ is supersymmetric (that is $\boldsymbol{A}$
is invariant under any permutation of it's indices) then $\boldsymbol{A}$
induces a single map. Furthermore, different constraints on the solution
eigenvectors $\boldsymbol{x}_{k}$ distinguishes the $E$-eigenvectors
from the $H$-eigenvectors and the $Z$-eigenvectors as introduced
and discussed in \cite{Qi_L2,Qi_L1}.

Our treatment considerably differs from the approaches described above
in the fact that our aim is to find a decomposition for a given tensor
$\boldsymbol{A}$ that provides a natural generalization for the concepts
of Hermitian and orthogonal matrices. Furthermore our approach is
not limited to supersymmetric tensors.

In connection with our investigations in the current work, we point
out another concepts from linear algebra for which the generalization
to tensor plays a significant role in complexity theory, that is the
notion of matrix rank. Indeed one may also find an extensive discussions
on the topic of tensor rank in \cite{Qi_L3,Grigoriev98,Hastad90,Raz2010,Hillar2009}.
The tensor rank problem is perhaps best described by the following
optimization problem. Given an $r$-tensor $\boldsymbol{A}=\left(a_{i_{1},\cdots,i_{r}}\right)$
we seek to solve the following problem which attempts to find an approximation
of $\boldsymbol{A}$ as a linear combination of rank one tensors.

\begin{equation}
\min_{\left(\otimes\boldsymbol{x}_{k}^{(t)}\right)_{1\le t\le r}\in\left(\bigotimes_{1\le t\le r}V_{t}\right)}\left\Vert \left(\sum_{1\le k\le l}\left(\lambda_{k}\right)^{r}\bigotimes_{1\le t\le r}\boldsymbol{x}_{k}^{(t)}\right)-\boldsymbol{A}\right\Vert 
\end{equation}

Our proposed tensor decomposition into lower order tensors relates
to the tensor rank problem but differs in the fact that the lower
order tensors arising from the spectral decomposition of $3$-tensors,
named eigen-matrices are not necessarily rank $1$ matrices.

\subsection{Existing Tensor Decomposition Framework}

Several approaches have been introduced for decomposing $k$-tensors
for $k\ge3$ in a way inspired by matrix SVD. SVD decomposes a matrix
$\boldsymbol{A}$ into $\boldsymbol{A}=\boldsymbol{U}\boldsymbol{\Sigma}\boldsymbol{V}^{T}$
and can be viewed as a decomposition of the matrix $\boldsymbol{A}$
into a summation of rank-1 matrices that can be written as 
\begin{equation}
\boldsymbol{A}=\sum_{i=1}^{r}\sigma_{i}\;\otimes(u_{i},v_{i})\label{eq:SVD}
\end{equation}
 where $r$ is the rank of $\boldsymbol{A}$, $u_{i},v_{i}$ are the
$i$-th columns of the orthogonal matrices $\boldsymbol{U}$ and $\boldsymbol{V}$,
and $\sigma_{i}$'s are the diagonal elements of $\boldsymbol{\Sigma}$,
i.e., the singular values. Here $\otimes(\cdot,\cdot)$ denotes the
outer product. The Canonical and Parallel factor decomposition (CANECOMP-PARAFAC,
also caller the CP model), independently introduced by \cite{Carroll70,Harshman70},
generalize the SVD by factorizing a tensor into a linear combination
of rank-1 tensors. That is given $\boldsymbol{A}\in\mathbb{R}^{n_{1}\times n_{2}\times n_{3}}$,
the goal is to find matrices $\boldsymbol{U}\in\mathbb{R}^{n_{1}\times n_{1}}$,
$\boldsymbol{V}\in\mathbb{R}^{n_{2}\times n_{2}}$ and $\boldsymbol{W}\in\mathbb{R}^{n_{3}\times n_{3}}$
such that 
\begin{equation}
\boldsymbol{A}=\sum_{i=1}^{r}\sigma_{i}\;\otimes(u_{i},v_{i},w_{i})\label{eq:rank1decomp}
\end{equation}
 where the expansion is in terms of the outer product of vectors $u_{i},v_{i},w_{i}$
are the i-th columns of $\boldsymbol{U}$, $\boldsymbol{V}$, and
$\boldsymbol{W}$, which yields rank-1 tensors. The rank of $\boldsymbol{A}$
is defined as the minimum $r$ required for such an expansion. Here
there are no assumption about the orthogonality of the column vectors
of $\boldsymbol{U}$, $\boldsymbol{V}$, and $\boldsymbol{W}$. The
CP decomposition have been show to be useful in several applications
where such orthogonality is not required. There are no known closed-form
solution to determine the rank $r$, or to find a lower rank approximation
as given directly by matrix SVD.

Tucker decomposition, introduced in\,\cite{Tucker66}, generalizes
over Eq \ref{eq:rank1decomp}, where an $\left(n_{1}\times n_{2}\times n_{3}\right)$
tensor $\boldsymbol{A}$ is decomposed into rank-1 tensor expansion
in the form
\begin{equation}
\boldsymbol{A}=\sum_{i=1}^{n_{1}}\sum_{j=1}^{n_{2}}\sum_{k=1}^{n_{3}}\sigma_{i,j,k}\;\otimes(u_{i},v_{j},w_{k})\label{eq:Tuckerdecomp}
\end{equation}
 where $u_{i}\in\mathbb{R}^{n_{1}}$, $v_{j}\in\mathbb{R}^{n_{2}}$,
and $w_{k}\in\mathbb{R}^{n_{3}}$. The coefficients $\sigma_{i,j,k}$
form a tensor that is called the core tensor $\boldsymbol{C}$. It
can be easily seen that if such core tensor is diagonal, i.e., $\sigma_{i,j,k}=0$
unless $i=j=k$, Tucker decomposition reduces to the CP decomposition
in Eq \ref{eq:rank1decomp}.

Orthogonality is not assumed in Tucker decomposition. Orthogonality
constraints can be added by requiring $u_{i},v_{j},w_{k}$ to be columns
of orthogonal matrices $\boldsymbol{U}$,$\boldsymbol{V}$, and $\boldsymbol{W}$.
Such decomposition was introduced in\,\cite{Lathauwer00JMAAa} and
was denoted by High Order Singular Value Decomposition (HOSVD). Tucker
decomposition can be written using the mode-$n$ tensor-matrix multiplication
defined in \cite{Lathauwer00JMAAa} as
\begin{equation}
\boldsymbol{A}=\boldsymbol{C}\times_{1}\boldsymbol{U}\times_{2}\boldsymbol{V}\times_{3}\boldsymbol{W}
\end{equation}
 where $\times_{n}$ is the mode-$n$ tensor-matrix multiplication.
Similar to Tucker decomposition, the core tensor of HOSVD is a dense
tensor. However, such a core tensor satisfies an \emph{all-orthogonality}
property between its slices across different dimensions as defined
in \cite{Lathauwer00JMAAa}.

HOSVD of a tensor can be computed by flattening the tensor into matrices
across different dimensions and using SVD on each matrix. Truncated
version of the expansion yields a lower rank approximation of a tensor
\cite{Lathauwer00JMAAb}. Several approaches have been introduced
for obtaining lower rank approximation by solving a least square problem,
\emph{e.g.} \cite{bb50202}. Recently an extension to Tucker decomposition
with non-negativity constraint was introduced with many successful
applications \cite{Shashua05NonNegative}.

All the above mentioned decompositions factorizes a high order tensor
as a summation of rank-1 tensors of the same dimension, which is inspired
by such an interpretation of matrix SVD as in Eq \ref{eq:SVD}. However,
none of these decomposition approaches can describe a tensor as a
product of tensors as would be expected from an SVD generalization.
The only known approach to us for decomposing a tensor to a product
of tensors was introduced in a technical report \cite{Kilmer08}.
This approach is based on the idea that a diagonalization of a circulant
matrix can be obtained by Discrete Fourier Transform (DFT). Given
a tensor, it is flattened then a block diagonal matrix is constructed
by DFT of the circulant matrix formed from the flattened tensor. Matrix
SVD is then used on each of the diagonal blocks. The inverse process
is then used to put back the resulting decompositions into tensors.
This approach results in a decomposition in the form $\boldsymbol{A}=\boldsymbol{U}\star\boldsymbol{S}\star\boldsymbol{V}^{T}$
where the product is defined as \cite{Kilmer08}
\[
\boldsymbol{A}\star\boldsymbol{B}=\mathtt{fold}(circ(\mathtt{unfold}(\boldsymbol{A},1)).\mathtt{unfold}(\boldsymbol{B},1),1)
\]
 However, such decomposition does not admit a representation of the
decomposition into an expansion in terms of rank-1 tensors. The product
is mainly defined by folding and unfolding the tensor into matrices.

From the above discussion we can highlight some fundamental limitations
of the known tensor decomposition frameworks. Existing tensor decomposition
frameworks are mainly expansions of a tensor as a linear combination
of rank-1 tensors, which are the outer products of vectors under certain
constraints (orthogonality, etc.) and do not provide a factorization
into product of tensors of the same dimensions. Tucker decomposition,
although a generalization of SVD, falls short of generalizing the
notion of the spectrum for high-order tensors. There is no connection
between the singular values and the spectrum of the corresponding
cubic Hermitian tensors. Unfortunately, no such relation is proposed
by the Tucker factorization. The Tucker decomposition does not suggest
at all how to generalize such objects as the trace and the determinant
of higher order tensors. In the appendix of this paper we show that
Tucker decomposition and HOSVD uses notion of matrix orthogonality.

\subsection{Applications of tensor decomposition\label{sub:Applications-of-tensor}}

The most widely used formulation for tensor decomposition is the orthogonal
version of Tucker decomposition (HOSVD) \cite{Lathauwer00JMAAa}.
HOSVD is a multilinear rank revealing procedure \cite{Lathauwer00JMAAa,Lathauwer00JMAAb}
and therefore, it has been widely used recently in many domains for
dimensionality reduction and to estimate signal subspaces of tensorial
data\,\cite{TensorReview}. In computer vision, HOSVD has been used
in\,\cite{Vasilescu02,Vasilescu03CVPR} for analysis of face images
with different sources of variability, \emph{e.g.} different people,
illumination, head poses, expressions, etc. It has been also used
in texture analysis, compression, motion analysis \cite{Vasilescu_CVPR01,Vasilescu_SIGGRAPH01},
posture estimation, gait biometric analysis, facial expression analysis
and synthesis, \emph{e.g.} \cite{Elgammal04CVPRa,LeeICCV07,Lee05AMFG,Lee05AVBPA},
and other useful applications \cite{TensorReview}. HOSVD decomposition
gives a natural way for dealing with images as matrices \cite{bb50202}.
The relation between HOSVD and independent component analysis ICA
was also demonstrated in \cite{LathauwerICA01} with applications
in communication, image processing, and others. Beyond vision and
image processing, HOSVD has also been used in data mining, web search,
\emph{e.g.} \cite{ICDM08,1106354,Sun05cubesvd:a}, and in DNA microarray
analysis \cite{TensorReview}.

\section{$3$-tensor algebra}

We propose a formulation for a general spectral theory for tensors
coined with consistent definitions from multilinear algebra. At the
core of the formulation is our proposed\emph{ spectral theory for
tensors} . In this section, the theory focuses on $3$-tensors algebra.
We shall discuss in the subsequent section the formulations of our
theory for $n$-tensor where $n$ is positive integer greater or equal
to $2$.

\subsection{Notation and Product definitions}

A $\left(m\times n\times p\right)$ \emph{3-tensor} $\boldsymbol{A}$
denotes a rectangular cuboid array of numbers having $m$ rows, $n$
columns, and $p$ depths. The entry $a_{i,j,k}$ occupies the position
where the $i^{th}$ row, the $j^{th}$ column, and the $k^{th}$ depth
meet. For many purposes it will suffice to write

\begin{equation}
\boldsymbol{A}\;:=\left(a_{i,j,k}\right)\;\left(1\le i\le m;\,1\le j\le n;\,1\le k\le p\right),
\end{equation}
We use the notation introduced above for matrices and vectors since
they will be considered special cases of $3$-tensors. Thereby, allowing
us to indicate matrices and vectors respectively as oriented slice
and fiber tensors. Therefore, $\left(m\times1\times1\right)$, $\left(1\times n\times1\right)$,
and $\left(1\times1\times p\right)$ tensors indicate vectors that
are respectively oriented vertically, horizontally and along the depth
direction furthermore they will be respectively denoted by $\boldsymbol{a}_{\centerdot,1,1}\;:=\left(a_{i,1,1}\right)_{\left\{ 1\le i\le m\right\} }$,
$\boldsymbol{a}_{1,\centerdot,1}\;:=\left(a_{1,j,1}\right)_{\left\{ 1\le j\le n\right\} }$,
$\boldsymbol{a}_{1,1,\centerdot}\;:=\left(a_{1,1,k}\right)_{\left\{ 1\le k\le p\right\} }$.
Similarly $\left(m\times n\times1\right)$, $\left(1\times n\times p\right)$,
and $\left(m\times1\times p\right)$ tensors indicate that the respective
martrices of dimensions $\left(m\times n\right)$, $\left(n\times p\right)$
and $\left(m\times p\right)$ can be respectively thought of as a
vertical, horizontal, or depth slice denoted respectively $\boldsymbol{a}_{\centerdot,\centerdot,1}\;:=\left(a_{i,j,1}\right)_{\left\{ 1\le i\le m,\:1\le j\le n\right\} }$,
$\boldsymbol{a}_{\centerdot,1,\centerdot}\;:=\left(a_{i,1,k}\right)_{\left\{ 1\le i\le m,\,1\le k\le p\right\} }$,
and $\boldsymbol{a}_{1,\centerdot,\centerdot}\;:=\left(a_{1,j,k}\right)_{\left\{ 1\le j\le n,\:1\le k\le p\right\} }$
.\\
\\
There are other definitions quite analogous to their matrix ($2$-tensors)
counterparts such as the definition of addition, Kronecker binary
product, and product of a tensor with a scalar, we shall skip such
definitions here.\\
\emph{Ternary product of tensors:} At the center of our proposed formulation
is the definition of the ternary product operation for $3$-tensors.
This definition, to the best of our knowledge has been first proposed
by P. Bhattacharya in \cite{Bhat} as a generalization of matrix multiplication.
Let $\boldsymbol{A}=\left(a_{i,j,k}\right)$ be a tensor of dimensions
$\left(m\times l\times p\right)$, $\boldsymbol{B}=\left(b_{i,j,k}\right)$
a tensor of dimensions $\left(m\times n\times l\right)$, and $\boldsymbol{C}=\left(c_{i,j,k}\right)$
a tensor of dimensions $\left(l\times n\times p\right)$; the ternary
product of $\boldsymbol{A}$, $\boldsymbol{B}$ and $\boldsymbol{C}$
results in a tensor $\boldsymbol{D}=\left(d_{i,j,k}\right)$ of dimensions
$(m\times n\times p)$ denoted 

\begin{equation}
\boldsymbol{D}=\circ\left(\boldsymbol{A},\boldsymbol{B},\boldsymbol{C}\right)
\end{equation}
and the product is expressed by :

\begin{equation}
d_{i,j,k}=\sum_{1\le t\le l}a_{i,t,k}\cdot b_{i,j,t}\cdot c_{t,j,k}\label{eq:ternaryproduct}
\end{equation}

\begin{figure}
\-\includegraphics[scale=0.3]{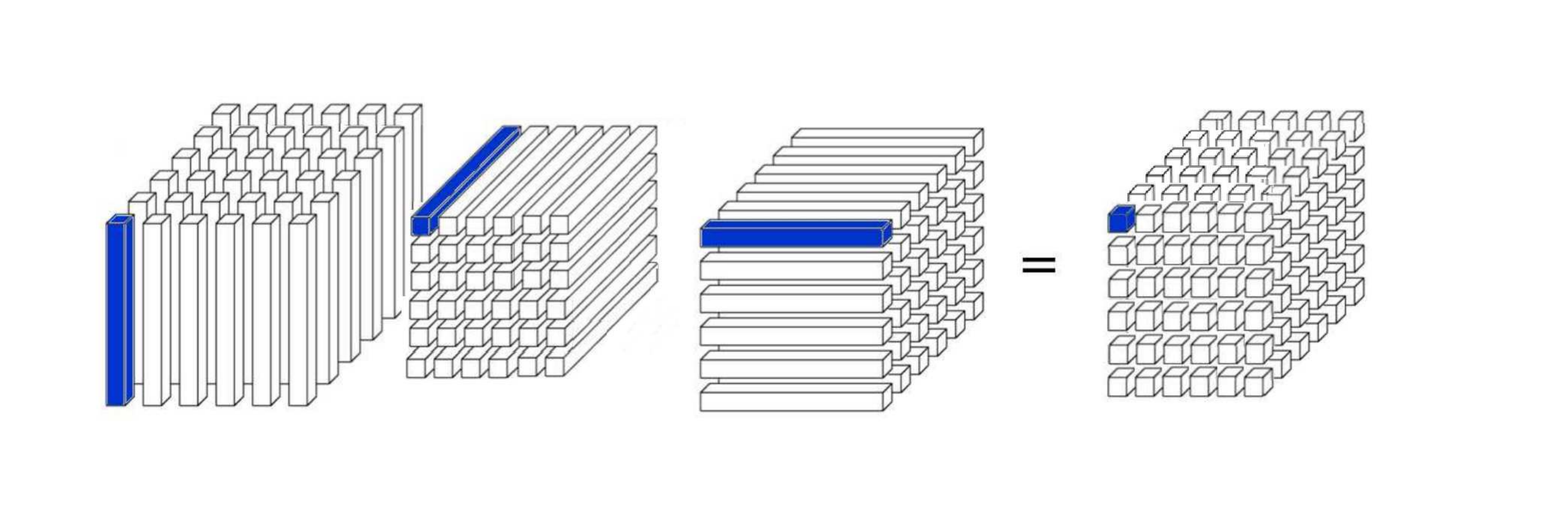}\caption{Tensor's ternary Product. }
\end{figure}

The specified dimensions of the tensors $\boldsymbol{A}$, $\boldsymbol{B}$
and $\boldsymbol{C}$ provide constraints for triplet of $3$-tensors
that can be multiplied using the preceding product definition. The
dimensions constraints are best illustrated by Fig.\,{[}2{]}. There
are several ways to generalize matrix product. We chose the previous
definition because the entries of the resulting tensor $\boldsymbol{D}=\circ\left(\boldsymbol{A},\boldsymbol{B},\boldsymbol{C}\right)$
relate to the general \emph{inner product} operator as depicted by
Fig.{[}1{]}. Therefore, the tensor product in Eq \ref{eq:ternaryproduct}
expresses the entries of $\boldsymbol{D}$ as inner products of the
triplet of horizontal, depth, and vertical vectors of $\boldsymbol{A}$,
$\boldsymbol{B}$ and $\boldsymbol{C}$ respectively as can be visualized
in Fig.\,{[}1{]}.

\begin{figure}
\includegraphics[scale=0.5]{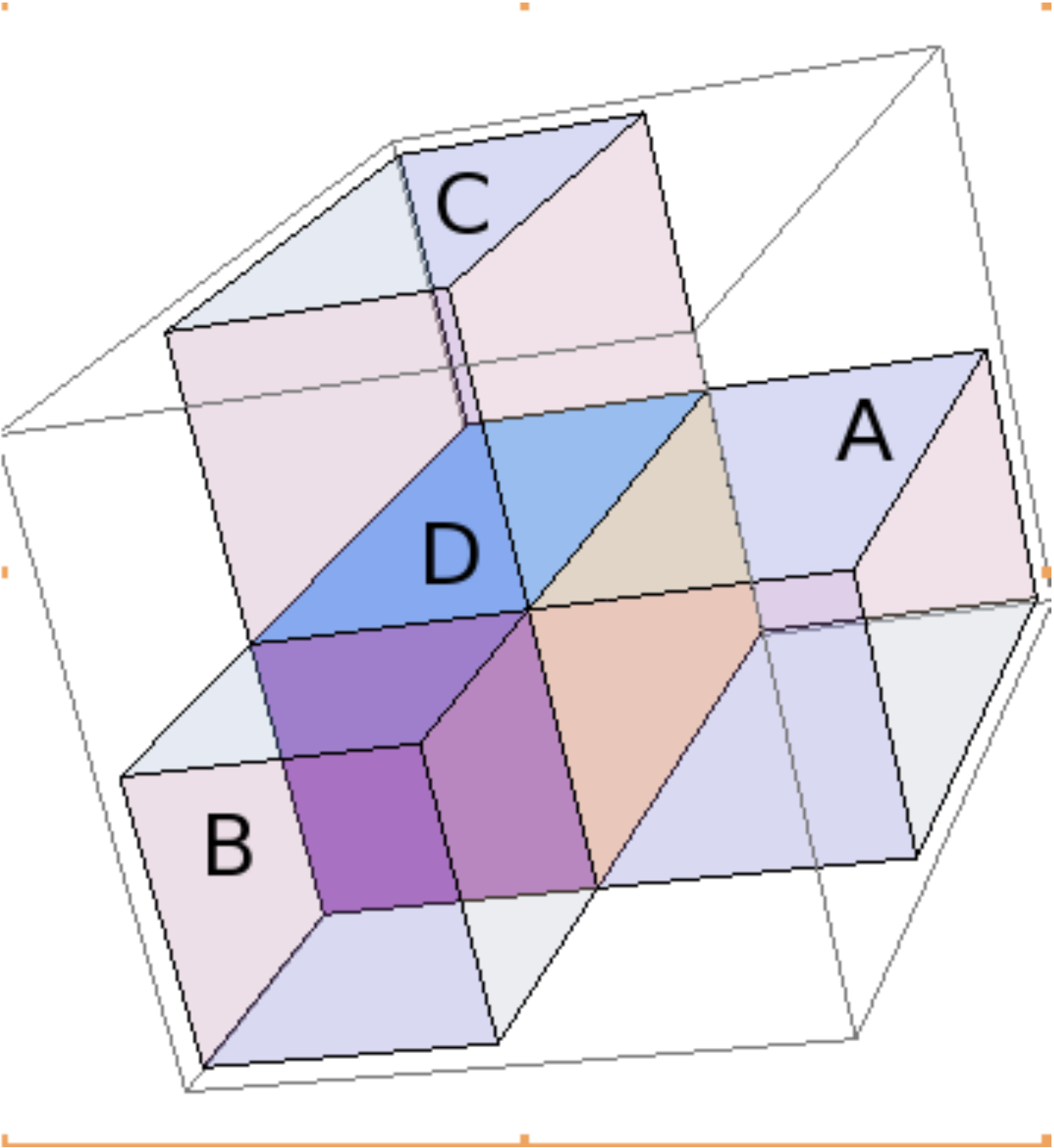}

\caption{Constraints on the dimensions of the tensors implied by the ternary
product definition.}
\end{figure}

We note that matrix product is a special instance of a tensor product
and we shall discuss subsequently products of $n$-tensor where $n$
is positive integer greater or equal to $2$. Furthermore the proposed
definition of the tensor multiplication suggests a generalization
of the binary vector outer product operator to a ternary operator
of slices. The ternary outer product is defined such that given tensors
$\boldsymbol{A}$ of dimensions $(m\times1\times p)$, $\boldsymbol{B}$
of dimensions $(m\times n\times1)$, and $\boldsymbol{C}$ of dimensions
$(1\times n\times p)$, their ternary outer product $\boldsymbol{D}$,
noted $\boldsymbol{D}=\otimes(\boldsymbol{A},\boldsymbol{B},\boldsymbol{C})$,
is an $(m\times n\times p)$ tensor defined by :

\begin{equation}
d_{i,j,k}=a_{i,1,k}\cdot b_{i,j,1}\cdot c_{1,j,k}.
\end{equation}
 Note that $\boldsymbol{A}$, $\boldsymbol{B},$ and $\boldsymbol{C}$
here are slices arising from oriented matrices. The above definition
generalizes the binary vector outer product operation to a ternary
matrix outer product operation defined by

\begin{equation}
\boldsymbol{D}=\otimes\left(\boldsymbol{a}_{\centerdot,1,\centerdot},\:\boldsymbol{b}_{\centerdot,\centerdot,1},\;\boldsymbol{c}_{1,\centerdot,\centerdot}\right)\;:=d_{i,j,k}=a_{i,1,k}\cdot b_{i,j,1}\cdot c_{1,j,k}.\label{eq:Outer-Product}
\end{equation}
Similarly to matrix multiplication, where the operation of multiplying
appropriate sized matrices can be viewed as a summation of outer product
of vectors, the product of appropriate sized triplet of tensors in
Eq \,\ref{eq:ternaryproduct} can be viewed as a summation of ternary
outer product of slices

\begin{equation}
\circ\left(\boldsymbol{A},\boldsymbol{B},\boldsymbol{C}\right)\equiv\sum_{1\le t\le l}\otimes\left(\boldsymbol{a}_{\centerdot,t,\centerdot},\:\boldsymbol{b}_{\centerdot,\centerdot,t},\;\boldsymbol{c}_{t,\centerdot,\centerdot}\right).
\end{equation}
\emph{}\\
\emph{Ternary dot product with a background tensor:} The ternary dot
product above can be further generalized by introducing the notion
of a background tensor as follows for $\boldsymbol{a}_{1,\centerdot,1}=\left(a_{1,i,1}\right)_{\left\{ 1\le i\le l\right\} }$,
$\boldsymbol{b}_{1,1,\centerdot}=\left(b_{1,1,j}\right)_{\left\{ 1\le j\le l\right\} }$
and $\boldsymbol{c}_{\centerdot,1,1}=\left(c_{k,1,1}\right)_{\left\{ 1\le k\le l\right\} }$

\begin{equation}
\left\langle \boldsymbol{a}_{1,\centerdot,1},\boldsymbol{b}_{1,1,\centerdot},\boldsymbol{c}_{\centerdot,1,1}\right\rangle _{\boldsymbol{T}}\;:=\sum_{1\le i\le l}\left(\sum_{1\le j\le l}\left(\sum_{1\le k\le l}a_{1,i,1}\cdot b_{1,1,j}^{\mathfrak{c}_{3}^{1}}\cdot c_{k,1,1}^{\mathfrak{c}_{3}^{2}}\cdot t_{i,j,k}\right)\right)\label{eq:Background-tensors}
\end{equation}
the preceding will be referred to as the \emph{triplet dot product}
operator with \emph{background tensor} $\boldsymbol{T}$. Background
tensors plays a role analogous to that of the metric tensor. The triplet
dot product with non trivial background tensor corresponds to a pure
trilinear form. Furthermore the outer product of $2$-tensors can
be generalized using the notion of background tensors to produce a
$3$-tensor $\boldsymbol{D}$ which result from a product of three
$2$-tensors namely $\boldsymbol{a}_{\centerdot,\centerdot,1}=\left(a_{m,i,1}\right)_{m,i}$,
$\boldsymbol{b}_{1,\centerdot,\centerdot}=\left(b_{1,n,j}\right)_{n,j}$
and $\boldsymbol{c}_{\centerdot,1,\centerdot}=\left(c_{k,1,p}\right)_{k,p}$
as follows,

\begin{equation}
d_{m,n,p}=\sum_{1\le i\le l}\left(\sum_{1\le j\le l}\left(\sum_{1\le k\le l}a_{m,i,1}\cdot b_{1,n,j}\cdot c_{k,1,p}\cdot t_{i,j,k}\right)\right).\label{eq:Matrices to Tensor}
\end{equation}
The preceding product expression is the one most commonly used as
a basis for tensor algebra in the literature as discussed in \cite{Hillar2009,Tucker66,LathauwerICA01,1106354}.\\
We may note that the original definition of the dot product for a
triplets of vectors corresponds to a setting where the background
tensor is the Kronecker delta $\boldsymbol{\Delta}=\left(\delta_{i,j,k}\right)$
that is $\boldsymbol{T}=\boldsymbol{\Delta}$ where $\boldsymbol{\Delta}$
denotes hereafter the Kronecker tensor and can be expressed in terms
of the Kronecker $2$-tensors as follows 

\begin{equation}
\delta_{i,j,k}=\delta_{i,j}\cdot\delta_{j,k}\cdot\delta_{k,i}
\end{equation}
equivalently $\boldsymbol{\Delta}=\left(\delta_{i,j,k}\right)$ can
be expressed in terms of the canonical basis $\left\{ \boldsymbol{e}_{i}:\;1\leq i\leq l\right\} $
in $l$-dimensional euclidean space described by:

\begin{equation}
\boldsymbol{\Delta}=\sum_{1\le k\le l}\left(\boldsymbol{e}_{k}\otimes\boldsymbol{e}_{k}\otimes\boldsymbol{e}_{k}\right),
\end{equation}

\begin{figure}
\center \includegraphics[scale=0.13]{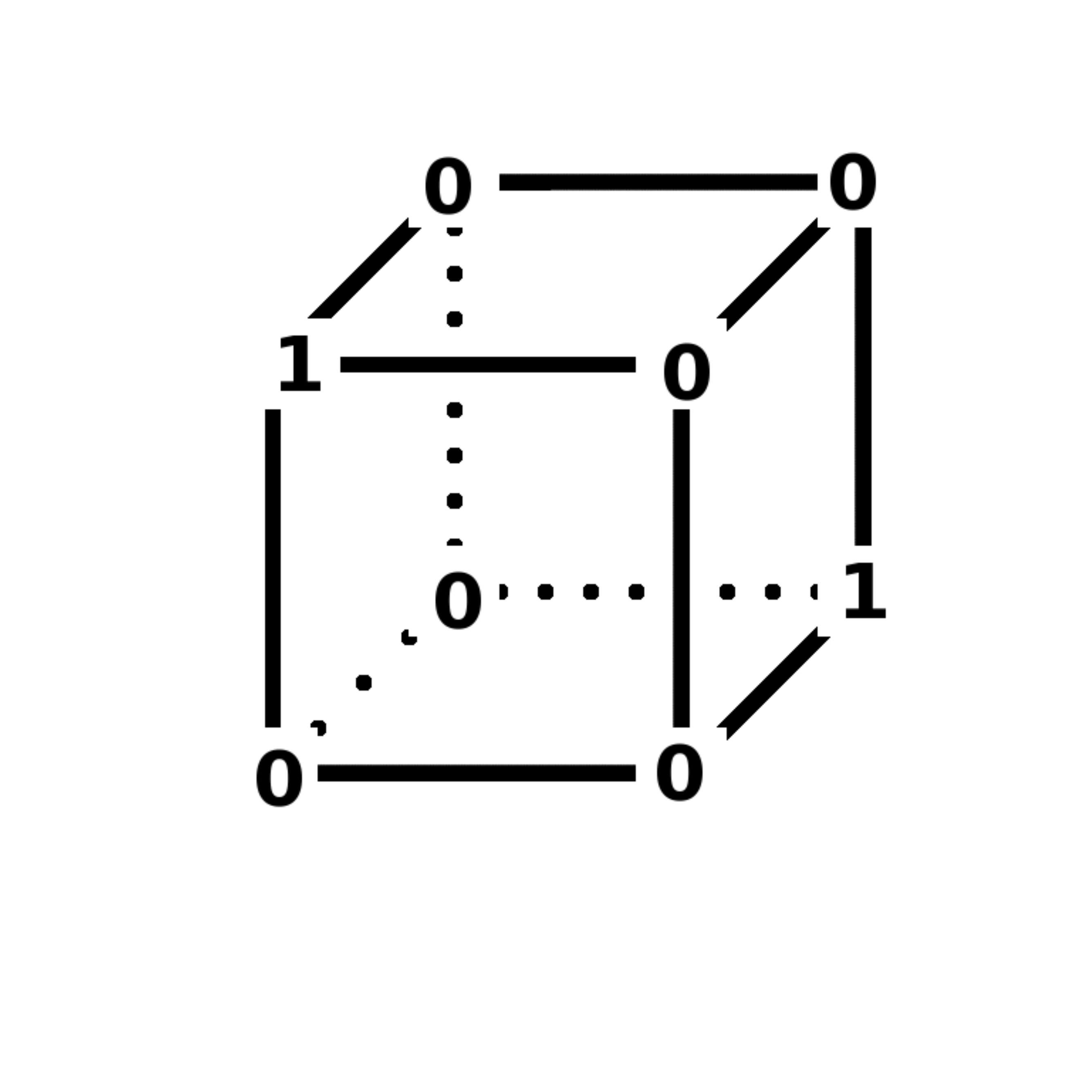}\caption{Kronecker $\left(2\times2\times2\right)$ tensor.}
\end{figure}

hence

\begin{equation}
\langle\boldsymbol{w},\,\boldsymbol{x},\,\boldsymbol{y}\rangle\equiv\langle\boldsymbol{w},\,\boldsymbol{x},\,\boldsymbol{y}\rangle_{\boldsymbol{\Delta}}.
\end{equation}

\subsubsection{Special Tensors and Special Operations}

In general it follows from the algebra described in the previous section
for $3$-tensors that:

\begin{equation}
\circ\left(\circ\left(\boldsymbol{A},\boldsymbol{B},\boldsymbol{C}\right),\boldsymbol{D},\boldsymbol{E}\right)\ne\circ\left(\boldsymbol{A},\circ\left(\boldsymbol{B},\boldsymbol{C},\boldsymbol{D}\right),\boldsymbol{E}\right)\ne\circ\left(\boldsymbol{A},\boldsymbol{B},\circ\left(\boldsymbol{C},\boldsymbol{D},\boldsymbol{E}\right)\right)\label{eq:Non_Associative}
\end{equation}
In some sense the preceding illustrates the fact that the product
operator is non associative over the set of tensors. However tensor
product is weakly distributive over tensor addition that is to say

\begin{equation}
\circ\left(\left[\boldsymbol{A}+\boldsymbol{B}\right],\,\boldsymbol{C},\,\boldsymbol{D}\right))=\circ\left(\boldsymbol{A},\boldsymbol{C},\boldsymbol{D}\right)+\circ\left(\boldsymbol{B},\boldsymbol{C},\boldsymbol{D}\right),
\end{equation}
however in general

\begin{equation}
\circ\left(\boldsymbol{A},\circ\left(\boldsymbol{B},\boldsymbol{C},\boldsymbol{D}\right),\boldsymbol{E}\right)+\circ\left(\boldsymbol{A},\circ\left(\boldsymbol{F},\boldsymbol{G},\boldsymbol{H}\right),\boldsymbol{E}\right)\ne\circ\left(\boldsymbol{A},\,\left(\circ\left(\boldsymbol{B},\boldsymbol{C},\boldsymbol{D}\right)+\circ\left(\boldsymbol{F},\boldsymbol{G},\boldsymbol{H}\right)\right)\,,\boldsymbol{E}\right).
\end{equation}
\emph{}\\
\emph{Transpose of a tensor:} Given a tensor $\boldsymbol{A}=\left(a_{u,v,w}\right)$
we define it's \emph{transpose} $\boldsymbol{A}^{T}$ and it's \emph{double
transpose }$\boldsymbol{A}^{T^{2}}$as follows:

\begin{equation}
\boldsymbol{A}^{T}=\left(a_{v,w,u}\right)\label{eq:Transpose}
\end{equation}

\begin{equation}
\boldsymbol{A}^{T^{2}}\equiv\left(\boldsymbol{A}^{T}\right)^{T}=\left(a_{w,u,v}\right).
\end{equation}
It immediately follows from the definition of the transpose that for
any tensor $\boldsymbol{A}$, $\boldsymbol{A}^{T^{3}}=\boldsymbol{A}.$
Incidentally the transpose operator corresponds to a cyclic permutation
of the indices of the entries of $\boldsymbol{A}$. Therefore we can
defined a inverse transpose $\boldsymbol{A}^{T^{-1}}=\boldsymbol{A}^{T^{2}}$,
generally we have

\begin{equation}
\boldsymbol{A}^{T^{q}}=\left(\boldsymbol{A}^{T^{q-1}}\right)^{T},
\end{equation}
furthermore, a tensor $\boldsymbol{A}$ is said to be symmetrical
if :

\begin{equation}
\boldsymbol{A}=\boldsymbol{A}^{T}=\boldsymbol{A}^{T^{2}}.
\end{equation}
As a result for a given arbitrary $3$-tensor $\boldsymbol{A}$, the
products $\boldsymbol{B}=\circ\left(\boldsymbol{A},\boldsymbol{A}^{T^{2}},\boldsymbol{A}^{T}\right)$,
$\boldsymbol{C}=\circ\left(\boldsymbol{A}^{T},\boldsymbol{A},\boldsymbol{A}^{T^{2}}\right)$
and $\boldsymbol{D}=\circ\left(\boldsymbol{A}^{T^{2}},\boldsymbol{A}^{T},\boldsymbol{A}\right)$
all result in symmetric tensors. It also follows from the definitions
of the transpose operation and the definition of ternary product operation
that:

\begin{equation}
\left[\circ\left(\boldsymbol{A},\boldsymbol{B},\boldsymbol{C}\right)\right]^{T}=\circ\left(\boldsymbol{B}^{T},\boldsymbol{C}^{T},\boldsymbol{A}^{T}\right)
\end{equation}
 and

\begin{equation}
\left[\circ\left(\boldsymbol{A},\boldsymbol{B},\boldsymbol{C}\right)\right]^{T^{2}}=\left[\circ\left(\boldsymbol{B}^{T},\boldsymbol{C}^{T},\boldsymbol{A}^{T}\right)\right]^{T}=\circ\left(\boldsymbol{C}^{T^{2}},\boldsymbol{A}^{T^{2}},\boldsymbol{B}^{T^{2}}\right).
\end{equation}
\\
\emph{Adjoint operator:} For $\boldsymbol{A}\in\mathbb{C}^{m\times n\times p}$
we introduce the analog of the \emph{adjoint} operator for $3$-tensors
in two steps. The first step consists in writing all the entries of
$\boldsymbol{A}$ in their complex polar form.

\begin{equation}
\boldsymbol{A}=\left(a_{u,v,w}=r_{u,v,w}\cdot\exp\left\{ i\cdot\theta_{u,v,w}\right\} \right)\;\left(1\le u\le m;\,1\le v\le n;\,1\le w\le p\right).
\end{equation}
 The final step expresses the adjoint of the tensor $\boldsymbol{A}$
noted $\boldsymbol{A}^{\dagger}$ as follows

\begin{equation}
\begin{cases}
\begin{array}{c}
\boldsymbol{A}^{\dagger}\equiv\left(\boldsymbol{A}^{\mathfrak{c}_{3}^{1}}\right)^{T}\::=\left(r_{v,w,u}\cdot\exp\left\{ i\,\exp\left\{ i\:\frac{2\pi}{3}\right\} \cdot\theta_{v,w,u}\right\} \right)\\
\boldsymbol{A}^{\dagger^{2}}\equiv\left(\boldsymbol{A}^{\mathfrak{c}_{3}^{2}}\right)^{T^{2}}\::=\left(r_{w,u,v}\cdot\exp\left\{ i\,\exp\left\{ i\:\frac{4\pi}{3}\right\} \cdot\theta_{w,u,v}\right\} \right)\\
\boldsymbol{A}^{\dagger^{3}}\equiv\left(\boldsymbol{A}^{\mathfrak{c}_{3}^{3}}\right)^{T^{3}}\::=\left(a_{u,v,w}=r_{u,v,w}\cdot\exp\left\{ i\cdot\theta_{u,v,w}\right\} \right)
\end{array} & .\end{cases}
\end{equation}
The adjoint operator introduced here allows us to generalize the notion
of Hermitian matrices or self adjoint matrices to tensors. A tensor
is Hermitian if the following identity holds

\begin{equation}
\boldsymbol{A}^{\dagger}=\boldsymbol{A}.
\end{equation}
Incidentally the products $\circ\left(\boldsymbol{A},\boldsymbol{A}^{\dagger^{2}},\boldsymbol{A}^{\dagger}\right)$,
$\circ\left(\boldsymbol{A}^{\dagger},\boldsymbol{A},\boldsymbol{A}^{\dagger^{2}}\right)$
and $\circ\left(\boldsymbol{A}^{\dagger^{2}},\boldsymbol{A}^{\dagger},\boldsymbol{A}\right)$
result in self adjoint tensors or Hermitian tensors.\\
\emph{Identity Tensor:} Let $\boldsymbol{1}_{(m\times n\times p)}$
denotes the tensor having all it's entries equal to one and of dimensions
$(m\times n\times p)$. Recalling that $\boldsymbol{\Delta}=\left(\delta_{i,j,k}\right)$
denotes the Kronecker 3-tensor, we define the \emph{identity tensors}
$\boldsymbol{I}$ to be :

\begin{equation}
\boldsymbol{I}=\circ\left(\boldsymbol{1}_{(l\times l\times l)},\boldsymbol{1}_{(l\times l\times l)},\boldsymbol{\Delta}\right)=\circ\left(\boldsymbol{1}_{(l\times l\times l)},\boldsymbol{1}_{(l\times l\times l)},\left(\sum_{1\le k\le l}\boldsymbol{e}_{k}\otimes\boldsymbol{e}_{k}\otimes\boldsymbol{e}_{k}\right)\right)
\end{equation}

\begin{equation}
\boldsymbol{I}\equiv\left(i_{m,n,p}=\left(\sum_{1\le k\le l}\delta_{k,n,p}\right)=\delta_{n,p}\right)
\end{equation}
 Furthermore we have :

\begin{equation}
\boldsymbol{I}^{T}=\circ\left(\boldsymbol{1}_{(l\times l\times l)},\boldsymbol{\Delta},\boldsymbol{1}_{(l\times l\times l)}\right)=\circ\left(\boldsymbol{1}_{(l\times l\times l)},\left(\sum_{1\le k\le l}\boldsymbol{e}_{k}\otimes\boldsymbol{e}_{k}\otimes\boldsymbol{e}_{k}\right),\boldsymbol{1}_{(l\times l\times l)}\right)
\end{equation}

\begin{equation}
\boldsymbol{I}^{T}\equiv\left(\left(\boldsymbol{I}^{T}\right)_{m,n,p}=\left(\sum_{1\le k\le l}\delta_{m,n,k}\right)=\delta_{m,n}\right)
\end{equation}

\begin{equation}
\boldsymbol{I}^{T^{2}}=\circ\left(\boldsymbol{\Delta},\boldsymbol{1}_{(l\times l\times l)},\boldsymbol{1}_{(l\times l\times l)}\right)=\circ\left(\left(\sum_{1\le k\le l}\boldsymbol{e}_{k}\otimes\boldsymbol{e}_{k}\otimes\boldsymbol{e}_{k}\right),\boldsymbol{1}_{(l\times l\times l)},\boldsymbol{1}_{(l\times l\times l)}\right)
\end{equation}

\begin{equation}
\boldsymbol{I}^{T^{2}}\equiv\left(\left(\boldsymbol{I}^{T^{2}}\right)_{m,n,p}=\left(\sum_{1\le k\le l}\delta_{m,k,p}\right)=\delta_{m,p}\right)
\end{equation}
 for all positive integer $l\ge2$ . The identity tensor plays a role
quite analogous to the role of the identity matrix since $\forall\,\boldsymbol{A}\in\mathbb{C}^{l\times l\times l}$
we have

\begin{equation}
\circ\left(\boldsymbol{I},\boldsymbol{A},\boldsymbol{I}^{T^{2}}\right)=\boldsymbol{A}.
\end{equation}
 \textbf{Proposition 1:} \emph{$\forall\boldsymbol{A}\quad\circ\left(\boldsymbol{X},\boldsymbol{A},\boldsymbol{X}^{T^{2}}\right)=\boldsymbol{A}\quad and\quad\boldsymbol{X=}\left(x_{m,n,p}\ge0\right)\Leftrightarrow\boldsymbol{X}=\boldsymbol{I}$}\\
We prove the preceding assertion in two steps, the first step consists
of showing that the $\boldsymbol{I}$ is indeed a solution to the
equation

\begin{equation}
\forall\boldsymbol{A}\quad\circ\left(\boldsymbol{X},\boldsymbol{A},\boldsymbol{X}^{T^{2}}\right)=\boldsymbol{A}
\end{equation}
Let $\boldsymbol{R}$ be the result of the product

\begin{equation}
\boldsymbol{R}=\left(r_{m,n,p}\right)=\circ\left(\boldsymbol{I},\boldsymbol{A},\boldsymbol{I}^{T^{2}}\right)
\end{equation}

\begin{equation}
r_{m,n,p}=\left(\sum_{1\le k\le l}i_{m,k,p}\cdot a_{m,n,k}\cdot\left(\boldsymbol{I}^{T^{2}}\right)_{k,n,p}\right)=\left(\sum_{1\le k\le l}\delta_{k,p}\cdot a_{m,n,k}\cdot\delta_{k,p}\right)
\end{equation}

\begin{equation}
r_{m,n,p}=\left(\sum_{1\le k\le l}\left(\delta_{k,p}\right)^{2}\cdot a_{m,n,k}\right)
\end{equation}
we note that

\begin{equation}
r_{m,n,k}=\begin{cases}
\begin{array}{c}
a_{m,n,k}\; if\: k=p\\
0\; otherwise
\end{array}\end{cases}
\end{equation}
hence

\begin{equation}
\boldsymbol{A}=\circ\left(\boldsymbol{I},\boldsymbol{A},\boldsymbol{I}^{T^{2}}\right).
\end{equation}
The last step consists in proving by contradiction that $\boldsymbol{I}$
is the unique solution with positive entries to the equation

\begin{equation}
\forall\boldsymbol{A}\quad\circ\left(\boldsymbol{X},\boldsymbol{A},\boldsymbol{X}^{T^{2}}\right)=\boldsymbol{A}
\end{equation}
Suppose there were some other solution $\boldsymbol{J}$ with positive
entry to the above equation, this would imply that

\begin{equation}
\circ\left(\boldsymbol{I},\boldsymbol{A},\boldsymbol{I}^{T^{2}}\right)-\quad\circ\left(\boldsymbol{J},\boldsymbol{A},\boldsymbol{J}^{T^{2}}\right)=0
\end{equation}

\begin{equation}
\Rightarrow\left(\sum_{1\le k\le l}i_{m,k,p}\cdot a_{m,n,k}\cdot\left(\boldsymbol{I}^{T^{2}}\right)_{k,n,p}\right)-\left(\sum_{1\le k\le l}j_{m,k,p}\cdot a_{m,n,k}\cdot\left(\boldsymbol{J}^{T^{2}}\right)_{k,n,p}\right)
\end{equation}

\begin{equation}
0=\sum_{1\le k\le l}a_{m,n,k}\cdot\left[\left(i_{m,k,p}\cdot\left(\boldsymbol{I}^{T^{2}}\right)_{k,n,p}\right)-\left(j_{m,k,p}\cdot\left(\boldsymbol{J}^{T^{2}}\right)_{k,n,p}\right)\right]
\end{equation}
Since this expression must be true for any choice of the values of
$a_{m,n,k}$ we deduce that it must be the case that

\begin{equation}
\left(\delta_{k,p}\right)^{2}-\left(j_{m,k,p}\cdot\left(\boldsymbol{J}^{T^{2}}\right)_{k,n,p}\right)=0
\end{equation}

\begin{equation}
\Rightarrow\left(j_{m,k,p}\cdot\left(\boldsymbol{J}^{T^{2}}\right)_{k,n,p}\right)=\delta_{k,p}
\end{equation}

\begin{equation}
j_{m,k,p}=\pm\delta_{k,p}
\end{equation}
the requirement that

\begin{equation}
j_{m,k,p}\ge0\Rightarrow j_{m,k,p}=\delta_{k,p}
\end{equation}
which results in the sought after contradiction $\square$.

\begin{figure}
\center \includegraphics[scale=0.13]{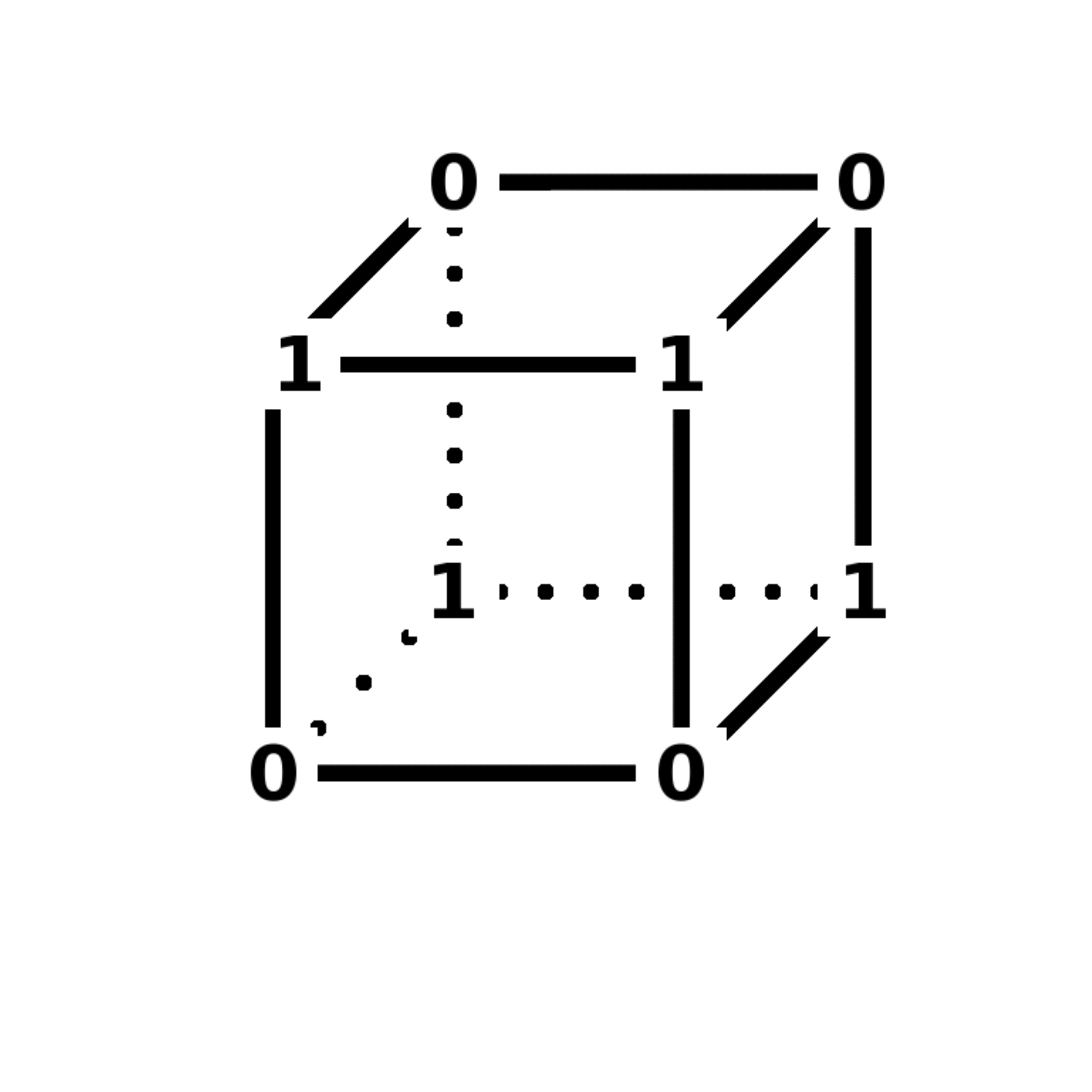} \includegraphics[scale=0.13]{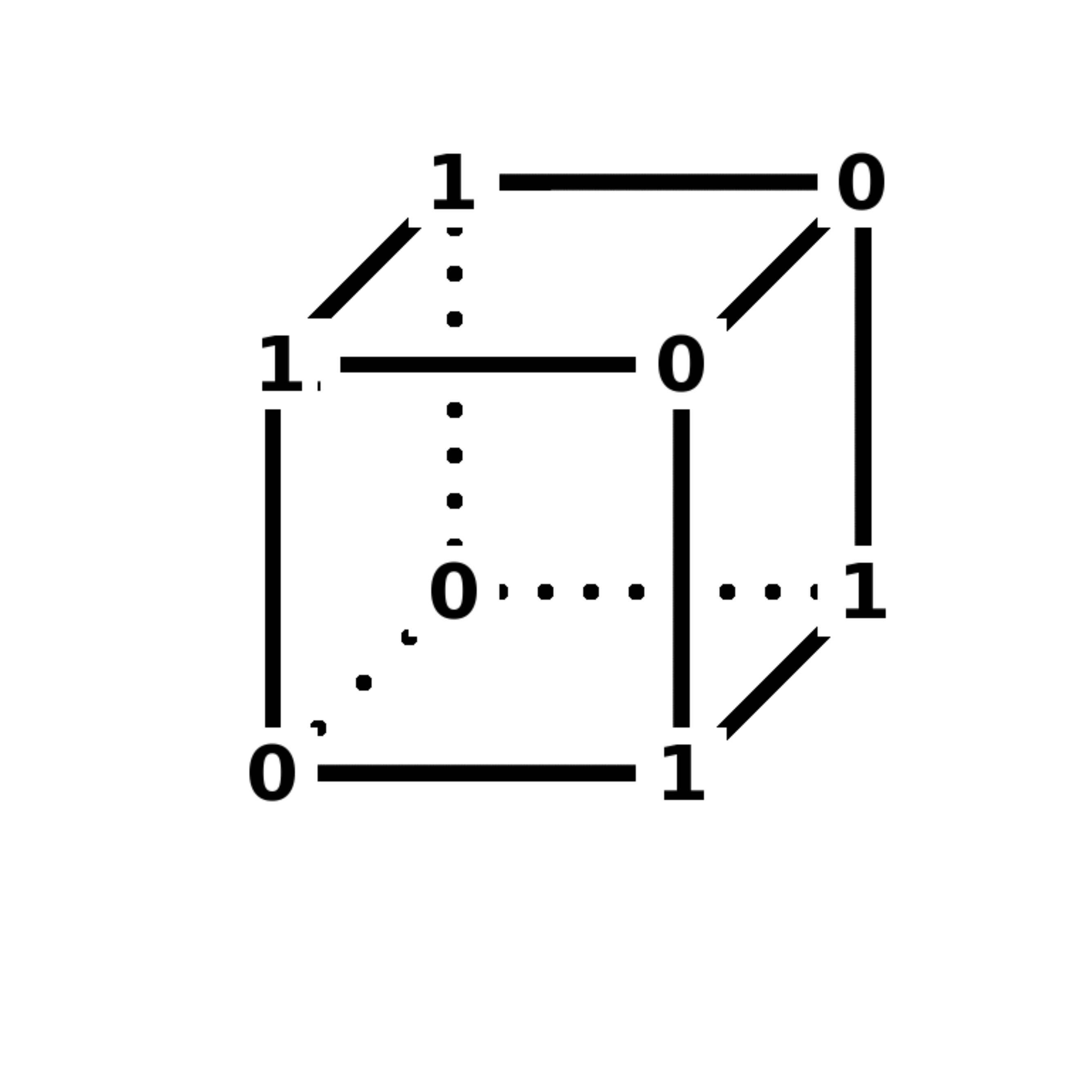}\includegraphics[scale=0.13]{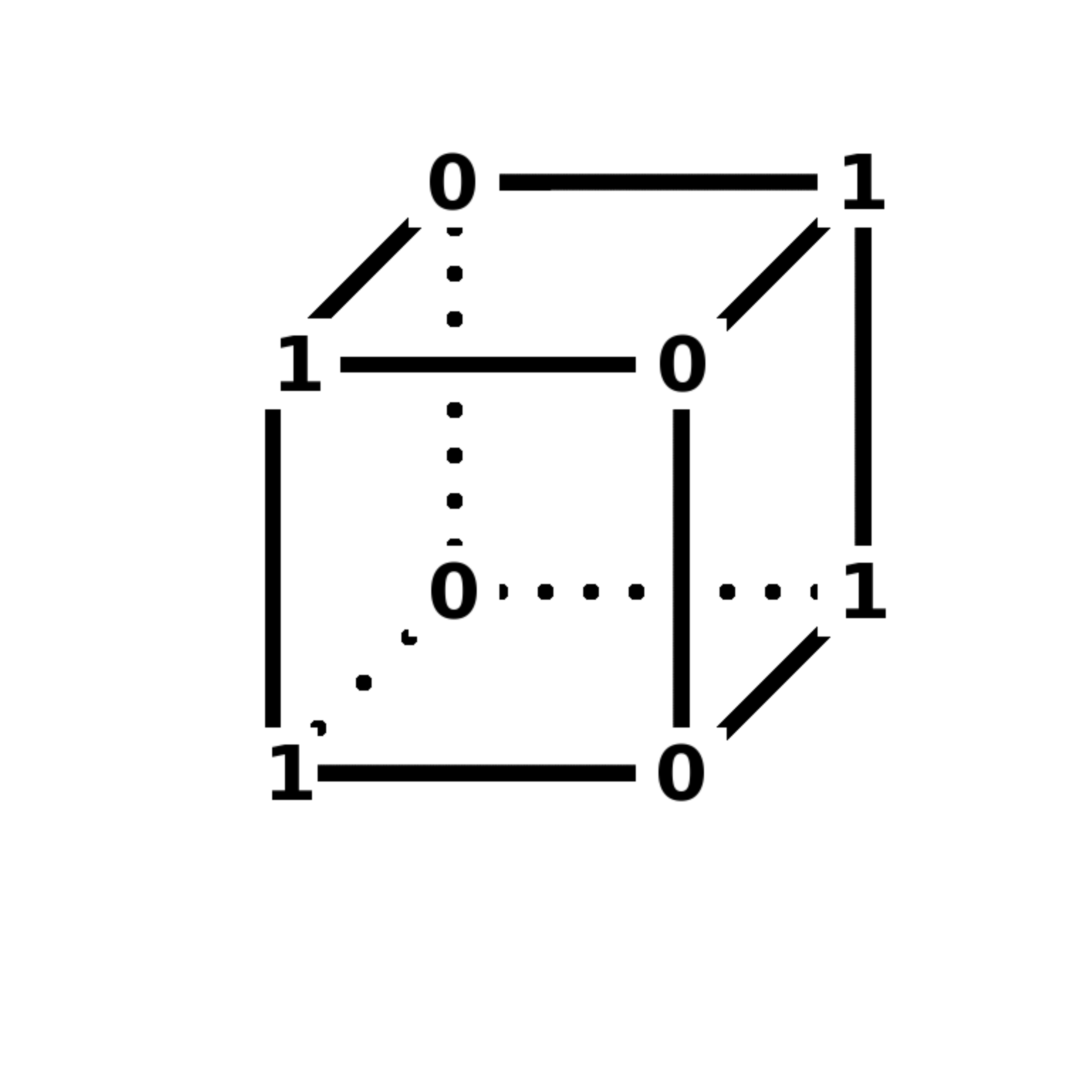}\caption{Tensor $\boldsymbol{I}$,$\boldsymbol{I}^{T}$and $\boldsymbol{I}^{T^{2}}$}
\end{figure}

\emph{Inverse:} By analogy to matrix inverse $\boldsymbol{A}^{-1}$
we recall that for a matrix $\boldsymbol{A}$, $\boldsymbol{A}^{-1}$
is its inverse if $\left(\boldsymbol{M}\boldsymbol{A}\right)\boldsymbol{A}^{-1}=\boldsymbol{M}$,
for any non zero matrix $\boldsymbol{M}$. We introduce here the notion
of inverse pairs for tensors. The ordered pair $\left(\boldsymbol{A}_{1},\boldsymbol{A}_{2}\right)$
and $\left(\boldsymbol{B}_{1},\boldsymbol{B}_{2}\right)$ are related
by inverse relationship if for any non zero $3$-tensor $\boldsymbol{M}$
with appropriated dimensions the following identity holds

\begin{equation}
\boldsymbol{M}=\circ\left(\boldsymbol{B}_{1}\circ\left(\boldsymbol{A}_{1},\boldsymbol{M},\boldsymbol{A}_{2}\right),\boldsymbol{B}_{2}\right).
\end{equation}
\emph{}\\
\emph{Permutation tensors:} Incidentally one may also discuss the
notion of \emph{permutation tensors} associated with any element $\sigma$
of the permutation group $S_{n}$.

\begin{equation}
\forall\;\sigma\in S_{n}\;\boldsymbol{P}_{\sigma}\equiv\circ\left(\boldsymbol{1}_{(n\times n\times n)},\boldsymbol{1}_{(n\times n\times n)},\left(\sum_{1\le k\le l}\boldsymbol{e}_{k}\otimes\boldsymbol{e}_{k}\otimes\boldsymbol{e}_{\sigma\left(k\right)}\right)\right)
\end{equation}

\begin{equation}
=\sum_{1\le k\le l}\circ\left(\boldsymbol{1}_{(n\times n\times n)},\boldsymbol{1}_{(n\times n\times n)},\left(\boldsymbol{e}_{k}\otimes\boldsymbol{e}_{k}\otimes\boldsymbol{e}_{\sigma\left(k\right)}\right)\right)
\end{equation}
The $3$-tensor $\boldsymbol{P}_{\sigma}$ perform the permutation
$\sigma$ on the depth slices of a $3$-tensor $\boldsymbol{A}$ through
the product $\circ\left(\boldsymbol{P}_{\sigma},\,\boldsymbol{A},\,\boldsymbol{P}_{\sigma}^{T^{2}}\right)$
, consequently the products $\circ\left(\boldsymbol{P}_{\sigma}^{T},\,\boldsymbol{P}_{\sigma}^{T^{2}},\boldsymbol{A}\right)$
and $\circ\left(\boldsymbol{A}\,,\boldsymbol{P}_{\sigma},\,\boldsymbol{P}_{\sigma}^{T}\right)$
perform the same permutation respectively on the row slices and the
column slices of $\boldsymbol{A}$.\\
 \\
 \textbf{Proposition 2:} Any permutation of the depth slices of $\boldsymbol{A}$
can be obtained by finite sequence of product of transposition, and
the sequence is of the form

\noindent 
\begin{equation}
\circ\left(\boldsymbol{P}\sigma_{n},\cdots,\circ\left(\boldsymbol{P}\sigma_{k},\cdots,\circ\left(\boldsymbol{P}\sigma_{1},\boldsymbol{A},\left(\boldsymbol{P}\sigma_{1}\right)^{T^{2}}\right)\cdots,\left(\boldsymbol{P}\sigma_{k}\right)^{T^{2}}\right),\cdots,\left(\boldsymbol{P}\sigma_{n}\right)^{T^{2}}\right).
\end{equation}
The preceding is easily verified using the definition above and the
permutation decomposition theorem \cite{Dummit}. Furthermore permutation
tensors suggest a generalization of bi-stochastic matrices to bi-stochastic
tensors through the Birkhoff-Von Neumann bi-stochastic matrix theorem.

\subsubsection{Orthogonality and scaling tensors}

From linear algebra we know that permutation matrices belong to both
the set of bi-stochastic matrices and to the set of orthogonal matrices.
We described above a approach for defining bi-stochastic $3$-tensors,
we shall address in this section the notion of orthogonality for $3$-tensors.
We recall from linear algebra that a matrix $\boldsymbol{Q}$ is said
to be orthogonal if

\begin{equation}
\boldsymbol{Q}^{\dagger}\cdot\boldsymbol{Q}=\boldsymbol{Q}\cdot\boldsymbol{Q}^{\dagger}=\boldsymbol{\Delta}.
\end{equation}
When we consider the corresponding equation for 3-tensors two distinct
interpretations arise. The first interpretation related to orthonormal
basis induced by the row or column vectors of the orthogonal matrix
$\boldsymbol{Q}$ that is :

\begin{equation}
\left\langle \boldsymbol{q}_{\centerdot,m},\boldsymbol{q}_{\centerdot,n}\right\rangle \equiv\left\langle \boldsymbol{q}_{m},\boldsymbol{q}_{n}\right\rangle =\left(\sum_{1\le k\le l}q_{k,m}\cdot q_{k,n}^{\mathfrak{c}_{2}^{1}}\right)=\delta_{m,n}
\end{equation}
The corresponding equation for a $3$-tensor $\boldsymbol{Q}=\left(q_{m,n,p}\right)$
of dimensions $\left(l\times l\times l\right)$ is given by:

\begin{equation}
\boldsymbol{\Delta}=\circ\left(\boldsymbol{Q},\boldsymbol{Q}^{\dagger^{2}},\boldsymbol{Q}^{\dagger}\right)\label{eq:Orthonormaility}
\end{equation}
or explicitly we can write:

\begin{equation}
\left\langle \boldsymbol{q}_{m,\centerdot,p},\boldsymbol{q}_{n,\centerdot,m},\boldsymbol{q}_{p,\centerdot,n}\right\rangle =\left(\sum_{1\le k\le l}q_{m,k,p}\cdot q_{n,k,m}^{\mathfrak{c}_{3}^{2}}\cdot q_{p,k,n}^{\mathfrak{c}_{3}^{1}}\right)=\delta_{m,n,p}.
\end{equation}
The second interpretation arises from the Kronecker invariance equation
expressed by:

\begin{equation}
\boldsymbol{\Delta}=\boldsymbol{Q}^{\dagger}\boldsymbol{\Delta}\boldsymbol{Q}=\left(\boldsymbol{Q}^{\dagger}\boldsymbol{\Delta}\boldsymbol{Q}\right)^{\dagger}.
\end{equation}
The corresponding Kronecker invariance equation for $3$-tensor is
given by :

\[
\boldsymbol{\Delta}=\circ\left(\circ\left(\boldsymbol{Q},\circ\left(\boldsymbol{Q}^{\dagger},\boldsymbol{Q}^{\dagger^{2}},\boldsymbol{\Delta}\right),\boldsymbol{Q}^{\dagger^{2}}\right),\boldsymbol{Q},\boldsymbol{Q}^{\dagger}\right)
\]

\begin{equation}
=\left[\circ\left(\circ\left(\boldsymbol{Q},\circ\left(\boldsymbol{Q}^{\dagger},\boldsymbol{Q}^{\dagger^{2}},\boldsymbol{\Delta}\right),\boldsymbol{Q}^{\dagger^{2}}\right),\boldsymbol{Q},\boldsymbol{Q}^{\dagger}\right)\right]^{\dagger}=\left[\circ\left(\circ\left(\boldsymbol{Q},\circ\left(\boldsymbol{Q}^{\dagger},\boldsymbol{Q}^{\dagger^{2}},\boldsymbol{\Delta}\right),\boldsymbol{Q}^{\dagger^{2}}\right),\boldsymbol{Q},\boldsymbol{Q}^{\dagger}\right)\right]^{\dagger^{2}}.\label{eq:Kronecker Invariance}
\end{equation}
While Kronecker invariance properly expresses a generalization of
the conjugation operation and the $3$-uniform hypergraph isomorphism
equation it does not follow from the first interpretation of orthogonality,
that is to say

\begin{equation}
\boldsymbol{\Delta}=\circ\left(\boldsymbol{Q},\boldsymbol{Q}^{\dagger^{2}},\boldsymbol{Q}^{\dagger}\right)\nRightarrow\circ\left(\circ\left(\boldsymbol{Q},\circ\left(\boldsymbol{Q}^{\dagger},\boldsymbol{Q}^{\dagger^{2}},\boldsymbol{\Delta}\right),\boldsymbol{Q}^{\dagger^{2}}\right),\boldsymbol{Q},\boldsymbol{Q}^{\dagger}\right)=\boldsymbol{\Delta}.
\end{equation}
We now discuss \emph{Scaling tensors. }The scaling\emph{ }tensor play
a role analogous to diagonal matrices in the fact that tensor multiplication
with scalling tensor results in a tensor whose vectors are scalled.
First we observe that the identity pairs of tensors should corespond
to special scaling tensors. The general family of diagonal tensors
are expressed by pairs of tensors $\boldsymbol{B}=\left(b_{m,n,p}\right)$
, $\boldsymbol{C}=\left(c_{m,n,p}\right)$ such that 

\begin{equation}
\boldsymbol{B}\equiv\left(b_{m,n,p}=\delta_{n,p}\cdot w_{p,m}\right)
\end{equation}

\begin{equation}
\boldsymbol{C}\equiv\left(c_{m,n,p}=\delta_{m,n}\cdot w_{m,p}\right)
\end{equation}
 The product $\boldsymbol{D}=\circ\left(\boldsymbol{A},\boldsymbol{B},\boldsymbol{C}\right)$
yields

\begin{equation}
d_{m,n,p}=\sum_{1\le k\le l}a_{m,k,p}\cdot\left(\delta_{n,k}\cdot w_{m,k}\right)\cdot\left(\delta_{k,n}\cdot w_{k,p}\right)
\end{equation}

\begin{equation}
\Rightarrow d_{m,n,p}=w_{m,n}\cdot a_{m,n,p}\cdot w_{n,p}
\end{equation}
 The expression above illustrates the fact that $w_{m,n}$ and $w_{n,p}$
scale the entry $a_{m,n,p}$ of the tensor $\boldsymbol{A}$, or equivalently
one may view the expression above as describing the non-uniform scaling
of the following vector $\left(a_{m,n,p}\right)_{1\le n\le l}$. The
vector scaling transform is expressed by

\begin{equation}
\left(a_{m,n,p}\right)_{1\le n\le l}\rightarrow\left(w_{m,n}\cdot a_{m,n,p}\cdot w_{n,p}\right)_{1\le n\le l}
\end{equation}
 Furthermore the scaling factors for a given vector may be viewed
as coming from the same vector of the scaling matrix $\boldsymbol{W}=\left(w_{m,n}\right)$
if the matrix $\boldsymbol{W}$ is symmetric. Finally we may emphasize
the analogy with diagonal matrices, which satisfy the following equation
independently of the value assigned to their non zero entries. For
a given $\boldsymbol{D}$, we solve for $\boldsymbol{C}$ such that

\begin{equation}
\left(\boldsymbol{D}\cdot\boldsymbol{C}\right)_{m,n}=d_{m,n}^{2}.
\end{equation}
We recall from matrix algebra that:

\begin{equation}
\boldsymbol{C}=\boldsymbol{D}
\end{equation}
 and furthermore

\begin{equation}
\boldsymbol{D}=\left(d_{m,n}=\delta_{m,n}\cdot w_{n}\right)
\end{equation}

\begin{equation}
\left(\boldsymbol{D}\cdot\boldsymbol{D}^{T}\right)_{m,n}=\begin{cases}
\begin{array}{c}
d_{m,n}^{2}\; if\; m=n\\
0\quad otherwise
\end{array}\end{cases}
\end{equation}
 By analogy we may define scaling tensors to be tensors satisfying
the following equation independently of the value of the nonzero tensors.

\begin{equation}
\left(a_{m,n,p}\right)^{3}=\sum_{1\le k\le l}a_{m,k,p}\cdot b_{m,n,k}\cdot c_{k,n,p}
\end{equation}
a possible solution is given by 

\begin{equation}
a_{m,n,p}=\delta_{m,p}\cdot w_{p,n}
\end{equation}

\begin{equation}
b_{m,n,p}=\delta_{n,p}\cdot w_{m,p}
\end{equation}

\begin{equation}
c_{m,n,p}=\delta_{m,n}\cdot w_{p,m}
\end{equation}
This is easily verified by computing the product

\begin{equation}
D=\circ\left(\boldsymbol{A},\boldsymbol{B},\boldsymbol{C}\right)\equiv d_{m,n,p}=\sum_{1\le k\le l}\left(\delta_{m,p}\cdot w_{p,k}\right)\cdot\left(\delta_{n,k}\cdot w_{m,k}\right)\cdot\left(\delta_{k,n}\cdot w_{p,k}\right)
\end{equation}

\begin{equation}
\Rightarrow d_{m,n,p}=\left(\delta_{m,p}\cdot w_{p,n}\right)\cdot\left(\delta_{n,n}\cdot w_{m,n}\right)\cdot\left(\delta_{n,n}\cdot w_{p,n}\right)
\end{equation}

\begin{equation}
\Rightarrow d_{m,n,p}=\left(\delta_{m,p}\cdot w_{p,n}\right)\cdot w_{m,n}\cdot w_{p,n}
\end{equation}

\begin{equation}
d_{m,n,p}=\begin{cases}
\begin{array}{c}
w_{m,n}^{3}\; if\; m=p\\
0\quad otherwise
\end{array}\end{cases}
\end{equation}
 Fig{[}4{]} provides an example of diagonal tensors. It so happens
that $\boldsymbol{A}$, $\boldsymbol{B}$, $\boldsymbol{C}$ discussed
above are related by transpose relation for third order tensors. This
fact considerably simplifies the formulation of the to diagonality
property common to both matrices and $3$-tensors. By analogy to matrices
we say for $3$-tensors that a tensor $\boldsymbol{D}=\left(d_{m,n,p}\right)$
is diagonal if independently of the value of the non zero entries
of $\boldsymbol{D}$ we have :

\[
\circ\left(\boldsymbol{D}^{T},\,\boldsymbol{D}^{T^{2}},\,\boldsymbol{D}\right)_{m,n,p}=d_{m,n,p}^{3}.
\]
\textbf{Proposition 3:} if a 3-tensor $\boldsymbol{D}$ can be expressed
in terms of a symmetric matrix $\boldsymbol{W}=\left(w_{m,n}=w_{n,m}\right)$
in the form $\boldsymbol{D}=\left(d_{m,n,p}=w_{m,n}\cdot\delta_{n,p}\right)$
then $\boldsymbol{D}$ is diagonal.\\
The proof of the proposition follows from the fact that :

\begin{equation}
\left(\boldsymbol{D}^{T}\right)_{m,n,p}=\left(w_{p,n}\cdot\delta_{n,m}\right)
\end{equation}

\begin{equation}
\left(\boldsymbol{D}^{T^{2}}\right)_{m,n,p}=\left(w_{n,p}\cdot\delta_{p,m}\right)
\end{equation}
 from which it follows that 

\begin{equation}
\circ\left(\boldsymbol{D}^{T},\,\boldsymbol{D}^{T^{2}},\,\boldsymbol{D}\right)_{m,n,p}=\left(w_{m,n}\right)^{3}\cdot\delta_{n,p}
\end{equation}

\section{Spectral Analysis of $3$-tensors}

\paragraph{Observations from the Eigen-Value/Vector equations.}

We briefly review well established properties of matrices and their
spectral decomposition, in order to emphasize how these properties
carry over to spectral decomposition of tensors. From the definition
of eigen-value/vector equation, we know that for a square hermitian
matrix $\boldsymbol{A}$, there must exist pairs of matrices $\boldsymbol{Q}$,
$\boldsymbol{R}$ and pairs of diagonal matrices $\boldsymbol{D}$,
$\boldsymbol{E}$ such that 

\begin{equation}
\begin{cases}
\begin{array}{ccc}
\boldsymbol{A} & = & \left(\boldsymbol{D}\boldsymbol{Q}\right)^{\dagger}\left(\boldsymbol{E}\boldsymbol{R}\right)\\
\boldsymbol{I} & = & \boldsymbol{Q}\boldsymbol{R}
\end{array} & ,\end{cases}\label{eq:matrix_decomp1}
\end{equation}
where the columns of $\boldsymbol{Q}^{\dagger}$ corresponds to the
left eigenvectors of $\boldsymbol{A}$, the rows of $\boldsymbol{R}$
corresponds to the right eigenvectors of $\boldsymbol{A}$ and the
entries of the diagonal matrix $\left(\boldsymbol{D}^{\dagger}\mathbf{E}\right)$
correspond to eigenvalues of $\boldsymbol{A}$. 

\begin{equation}
a_{m,n}=\sum_{1\le k\le l}\left(\mu_{k}\, q_{k,m}\right)^{\mathfrak{c}_{2}^{1}}\left(\nu_{k}\, r_{k,n}\right).\label{eq:matrix_decomp2}
\end{equation}
Let $f_{m,n}(k)=q_{k,m}^{\mathfrak{c}_{3}^{2}}\cdot r_{k,n}$, i.e.,
the entries of the matrix resulting from the outer product of the
$k$-th left eigenvector with the $k$-th right eigenvector, incidentally
the spectral decomposition yields the following expansion which is
crucial to the principal component analysis scheme.

\begin{equation}
a_{m,n}=\sum_{1\le k\le l}\left(\mu_{k}^{\mathfrak{c}_{2}^{1}}\cdot\nu_{k}\right)\quad f_{m,n}(k)
\end{equation}
The preceding amounts to assert that the spectral decomposition offers
for every entry of the $2$-tensor $\boldsymbol{A}$ a positional
encoding in a basis formed by the eigenvalues of the matrix. Assuming
that the eigenvalues are sorted in decreasing order, the preceding
expression suggest an approximation scheme for the entries of $\boldsymbol{A}$
and, therefore, an approximation scheme for the $2$-tensor $\boldsymbol{A}$
itself.

\paragraph*{Definition}

The spectrum of an $n$-tensor corresponds to the collection of lower
order tensors the entry of which are solutions to the characteristic
system of equations.

\paragraph{Spectrum of Hermitian tensors}

The aim of this section is to rigorously characterize the spectrum
of a symmetric tensor of dimensions $(l\times l\times l)$. Fig.\,{[}5{]}
depicts the product and the slice that will subsequently also be referred
to as eigen-matrices.

\begin{figure}
\includegraphics[scale=0.2]{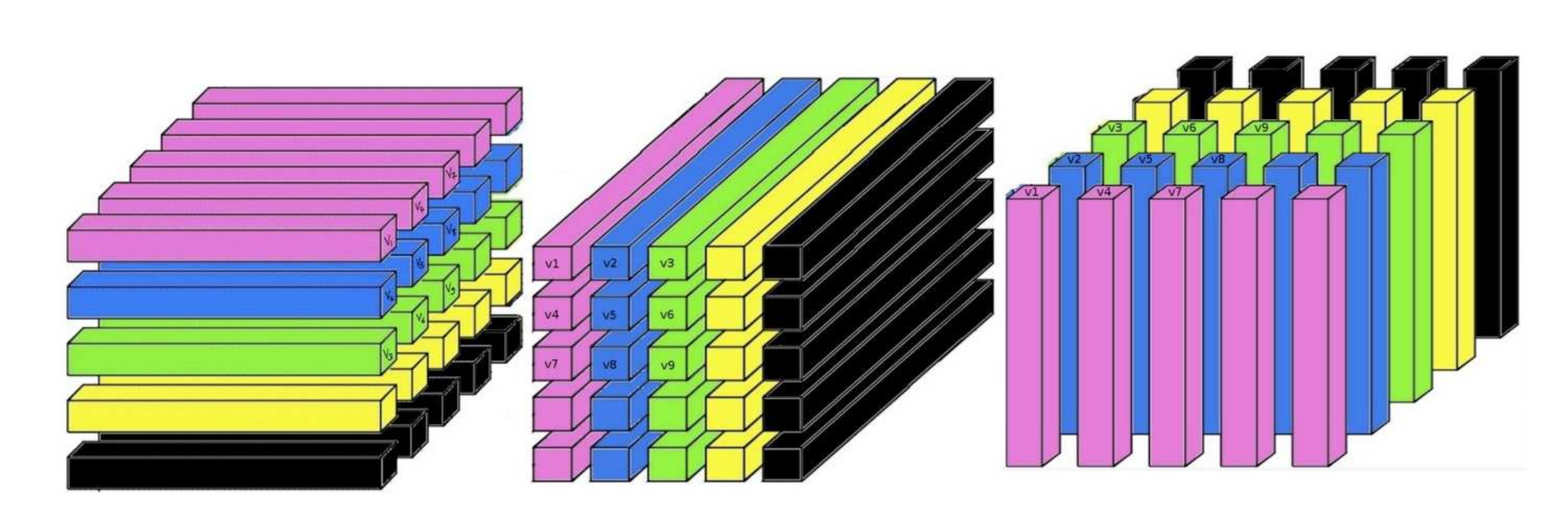}\caption{Orthogonal slices of an orthogonal tensor}
\end{figure}

We may state the \emph{spectral theorem} as follows\\

\noindent \textbf{Theorem 1:}\emph{ (Spectral Theorem for $3$-Tensors):}
For an arbitrary hermitian non zero $3$-tensor $\boldsymbol{A}$
with $\left\Vert \boldsymbol{A}\right\Vert _{\ell_{3}}^{3}\ne1$ there
exist a factorization of the form:

\begin{equation}
\begin{cases}
\begin{array}{ccc}
\boldsymbol{A} & = & \circ\left(\circ\left(\boldsymbol{Q},\boldsymbol{D},\boldsymbol{D}^{T}\right),\left[\circ\left(\boldsymbol{R},\boldsymbol{E},\boldsymbol{E}^{T}\right)\right]^{\dagger^{2}},\left[\circ\left(\boldsymbol{S},\boldsymbol{F},\boldsymbol{F}^{T}\right)\right]^{\dagger}\right)\\
\boldsymbol{\Delta} & = & \circ\left(\boldsymbol{Q},\:\boldsymbol{R}^{\dagger^{2}},\:\boldsymbol{S}^{\dagger}\right)
\end{array}\end{cases}\label{eq:3tensor_factorization}
\end{equation}
where $\boldsymbol{D}$, $\boldsymbol{E}$, $\boldsymbol{F}$ denote
scaling tensors. For convenience we introduce the following notation
for scaled tensors 

\begin{equation}
\begin{cases}
\begin{array}{ccc}
\widetilde{\boldsymbol{Q}} & = & \circ\left(\boldsymbol{Q},\boldsymbol{D},\boldsymbol{D}^{T}\right)\\
\widetilde{\boldsymbol{R}} & = & \circ\left(\boldsymbol{R},\boldsymbol{E},\boldsymbol{E}^{T}\right)\\
\widetilde{\boldsymbol{S}} & = & \circ\left(\boldsymbol{S},\boldsymbol{F},\boldsymbol{F}^{T}\right)
\end{array}\end{cases}\label{eq:Orthogonal Scaling}
\end{equation}
and simply expresses the tensor decomposition of $\boldsymbol{A}$
as:

\begin{equation}
\boldsymbol{A}=\circ\left(\widetilde{\boldsymbol{Q}},\widetilde{\boldsymbol{R}}^{\dagger^{2}},\widetilde{\boldsymbol{S}}^{\dagger}\right)\label{eq:Tensor factorization scaled}
\end{equation}

\subsection{Proof of the Spectral Theorem}

In what follows the polynomial ideal generated by the set of polynomials
$\left\{ f_{k}\right\} _{1\le k\le N}$ is noted $\left\rangle f_{k}\right\langle _{1\le k\le N}$.
We first emphasize the similarity between the spectral theorem for
tensors and matrices, by providing an alternative proof of a weaker
form of the spectral theorem for hermitian matrices with Forbenius
norm different from $1$. Finally we extend the proof technic to $3$-tensors
and subsequently to $n$-tensors.

\subsubsection*{Proof of the weak form of the spectral theorem for matrices}

Our aim is to prove that the spectral decomposition exists for an
arbitrary matrix $\boldsymbol{A}$ with forbenius norm different $1$.
For this we consider the ideals induced by the characteristic system
of equations for matrices. The spectral decomposition of $\boldsymbol{A}$
refers to the decomposition:

\begin{equation}
\begin{cases}
\begin{array}{ccc}
\boldsymbol{A} & = & \left(\boldsymbol{D}\boldsymbol{Q}\right)^{\dagger}\left(\boldsymbol{E}\boldsymbol{R}\right)\\
\boldsymbol{I} & = & \boldsymbol{Q}\boldsymbol{R}
\end{array} & ,\end{cases}
\end{equation}
the spectral decomposition equation above provides us with polynomial
system of equations in the form

\begin{equation}
\begin{cases}
\begin{array}{ccc}
a_{m,n} & = & \sum_{1\le k\le l}\left(\mu_{k}\, q_{k,m}\right)^{\mathfrak{c}_{2}^{1}}\left(\nu_{k}\, r_{k,n}\right)\\
\delta_{m,n} & = & \sum_{1\le k\le l}q_{k,m}^{\mathfrak{c}_{2}^{1}}\cdot r_{k,n}
\end{array} & 1\le m\le n\le l\end{cases}\label{eq:Characteristic_System}
\end{equation}
conveniently rewritten as

\begin{equation}
\begin{cases}
\begin{array}{c}
\left\langle \boldsymbol{D}\cdot\boldsymbol{q}_{m},\boldsymbol{E}\cdot\boldsymbol{r}_{n}\right\rangle =a_{m,n}\\
\left\langle \boldsymbol{q}_{m},\boldsymbol{r}_{n}\right\rangle =\delta_{m,n}
\end{array} & 1\le m\le n\le l\end{cases}.
\end{equation}
The ideal being considered is :

\begin{equation}
\mathcal{I}=\left\rangle \;\left\langle \boldsymbol{D}\cdot\boldsymbol{q}_{m},\boldsymbol{E}\cdot\boldsymbol{r}_{n}\right\rangle -a_{m,n},\:\left\langle \boldsymbol{q}_{m},\boldsymbol{r}_{n}\right\rangle -\delta_{m,n}\:\right\langle _{1\le m\le n\le l}\subseteq\mathbb{C}\left[\left\{ \mu_{k},\nu_{k}\:;\:\boldsymbol{q}_{k},\boldsymbol{r}_{k}\right\} _{1\le k\le l}\right].
\end{equation}
where the variables are the entries of the pairs of matrices $\boldsymbol{Q}$,
$\boldsymbol{R}$ and 
\[
\begin{array}{c}
\boldsymbol{D}=\left(d_{m,n}=\delta_{m,n}\cdot\mu_{m}\right)\\
\boldsymbol{E}=\left(e_{m,n}=\delta_{m,n}\cdot\nu_{m}\right)
\end{array}
\]
\textbf{Weak Spectral Theorem }\emph{(for $2$-tensors)}\textbf{:}
For an arbitrary non zero hermitian $2$-tensor $\boldsymbol{A}$
with $\left\Vert \boldsymbol{A}\right\Vert _{\ell_{2}}\ne1$ the spectral
system of polynomial equations :
\begin{equation}
\begin{cases}
\begin{array}{ccc}
\left\langle \boldsymbol{D}\cdot\boldsymbol{q}_{m},\boldsymbol{E}\cdot\boldsymbol{r}_{n}\right\rangle  & = & a_{m,n}\\
\left\langle \boldsymbol{q}_{m},\boldsymbol{r}_{n}\right\rangle  & = & \delta_{m,n}
\end{array} & 1\le m\le n\le l\end{cases}
\end{equation}
admits a solution.

\subsubsection*{\emph{Proof : }}

We prove this theorem by exhibiting a polynomial $p\left(\boldsymbol{D},\boldsymbol{E},\:\boldsymbol{q}_{1},\boldsymbol{r}_{1},\cdots,\boldsymbol{q}_{l},\boldsymbol{r}_{l}\right)$
which does not belong to the following ideal 
\[
\mathcal{I}=\left\rangle \;\left\langle \boldsymbol{D}\cdot\boldsymbol{q}_{m},\boldsymbol{E}\cdot\boldsymbol{r}_{n}\right\rangle -a_{m,n};\:\left\langle \boldsymbol{q}_{m},\boldsymbol{r}_{n}\right\rangle -\delta_{m,n}\;\right\langle _{1\le m,n\le l}.
\]
Consider the polynomial
\begin{equation}
p\left(\boldsymbol{D},\boldsymbol{E},\:\boldsymbol{q}_{1},\boldsymbol{r}_{1},\cdots,\boldsymbol{q}_{l},\boldsymbol{r}_{l}\right)\::=\left(\sum_{1\le m,n\le l}\left|\left\langle \boldsymbol{D}\cdot\boldsymbol{q}_{m},\boldsymbol{E}\cdot\boldsymbol{r}_{n}\right\rangle \right|^{2}\right)^{2}-\left\Vert \boldsymbol{A}\right\Vert _{\ell_{2}}^{2}.
\end{equation}
We claim that 
\begin{equation}
p\left(\boldsymbol{D},\boldsymbol{E},\:\boldsymbol{q}_{1},\boldsymbol{r}_{1},\cdots,\boldsymbol{q}_{l},\boldsymbol{r}_{l}\right)\notin\mathcal{I}
\end{equation}
since
\begin{equation}
p\left(\boldsymbol{D},\boldsymbol{E},\:\boldsymbol{q}_{1},\boldsymbol{r}_{1},\cdots,\boldsymbol{q}_{l},\boldsymbol{r}_{l}\right)\in\mathcal{I}\Rightarrow\left\Vert \boldsymbol{A}\right\Vert _{\ell_{2}}^{2^{2}}=\left\Vert \boldsymbol{A}\right\Vert _{\ell_{2}}^{2}
\end{equation}
which contradicts to the assumption that $\left\Vert \boldsymbol{A}\right\Vert _{\ell_{2}}^{2}\ne1$.
Hence we conclude that 
\begin{equation}
\left\Vert \boldsymbol{A}\right\Vert _{\ell_{2}}^{2}\ne1\Rightarrow p\left(\boldsymbol{D},\boldsymbol{E},\:\boldsymbol{q}_{1},\boldsymbol{r}_{1},\cdots,\boldsymbol{q}_{l},\boldsymbol{r}_{l}\right)\notin\mathcal{I}
\end{equation}
which completes the proof. $\square$ \\
In the proof above hermicity played a crucial role in that it ensures
that the eigenvalues are not all zeros since for non zero hermitian
$2$-tensor $\boldsymbol{A}$ 
\begin{equation}
\left\Vert \mathbf{A}\right\Vert _{\ell_{2}}^{2}=\mbox{Tr}\left\{ \mathbf{A}\cdot\mathbf{A}\right\} >0
\end{equation}

\subsubsection*{Proof of the Spectral Theorem for $3$-tensors}

We procede to derive the existence of spectral decomposition for $3$-tensors
using the proof thechnic discussed above

\begin{equation}
\begin{cases}
\begin{array}{ccc}
\boldsymbol{A} & = & \circ\left(\circ\left(\boldsymbol{Q},\boldsymbol{D},\boldsymbol{D}^{T}\right),\left[\circ\left(\boldsymbol{R},\boldsymbol{E},\boldsymbol{E}^{T}\right)\right]^{\dagger^{2}},\left[\circ\left(\boldsymbol{S},\boldsymbol{F},\boldsymbol{F}^{T}\right)\right]^{\dagger}\right)\\
\boldsymbol{\Delta} & = & \circ\left(\boldsymbol{Q},\:\boldsymbol{R}^{\dagger^{2}},\:\boldsymbol{S}^{\dagger}\right)
\end{array}\end{cases}
\end{equation}
equivalently written as 
\begin{equation}
\begin{cases}
\begin{array}{c}
a_{m,n,p}=\sum_{k=1}^{l}\left(\mu_{m,k}\cdot q_{m,k,p}\cdot\mu_{k,p}\right)\cdot\left(\nu_{n,k}\cdot r_{n,k,m}\cdot\nu_{k,m}\right)^{\mathfrak{c}_{3}^{2}}\cdot\left(\xi_{p,k}\cdot s_{p,k,n}\cdot\xi_{k,n}\right)^{\mathfrak{c}_{3}^{1}}\\
\delta_{m,n,p}=\sum_{k=1}^{l}q_{m,k,p}\cdot r_{n,k,m}^{\mathfrak{c}_{3}^{2}}\cdot s_{p,k,n}^{\mathfrak{c}_{3}^{1}}
\end{array}\end{cases}.
\end{equation}
The variables in the polynomial system of equations are the entries
of the $3$-tensor $\boldsymbol{Q}$, $\boldsymbol{R}$, $\boldsymbol{S}$
and the entries of the scaling tensors $\boldsymbol{D}$, $\boldsymbol{E}$,
$\boldsymbol{F}$ . \\
It is somewhat insightfull to express the system of equations in a
similar form to that of matrix spectral system of equations using
inner product moperators : 
\begin{equation}
\begin{cases}
\begin{array}{c}
\left\langle \boldsymbol{D}_{m,p}\cdot\boldsymbol{q}_{m,\centerdot,p},\:\boldsymbol{E}_{n,m}\cdot\boldsymbol{r}_{n,\centerdot,m},\:\boldsymbol{F}_{p,n}\cdot\boldsymbol{s}_{p,\centerdot,n}\right\rangle =a_{m,n,p}\\
\left\langle \boldsymbol{q}_{m,\centerdot,p},\,\boldsymbol{r}_{n,\centerdot,m},\,\boldsymbol{s}_{p,\centerdot,n}\right\rangle =\delta_{m,n,p}
\end{array}\end{cases}
\end{equation}
where $\boldsymbol{D}_{u,v}$ is a diagonal matrix whose entries are
specified by 
\begin{equation}
\begin{cases}
\begin{array}{c}
\boldsymbol{D}_{u,v}=\left(d_{i,j}^{u,v}=\delta_{i,j}\mu_{u,i}\mu_{j,v}\right)\\
\boldsymbol{E}_{u,v}=\left(e_{i,j}^{u,v}=\delta_{i,j}\nu_{u,i}\nu_{j,v}\right)\\
\boldsymbol{F}_{u,v}=\left(f_{i,j}^{u,v}=\delta_{i,j}\xi_{u,i}\xi_{j,v}\right)
\end{array}\end{cases}
\end{equation}
The characteristic system of equations yields the ideal $\mathcal{I}$
defined by

\begin{equation}
\mathcal{I}=\left\rangle \:\left\langle \boldsymbol{D}_{m,p}\cdot\boldsymbol{q}_{m,\centerdot,p},\:\boldsymbol{E}_{n,m}\cdot\boldsymbol{r}_{n,\centerdot,m},\:\boldsymbol{F}_{p,n}\cdot\boldsymbol{s}_{p,\centerdot,n}\right\rangle -a_{m,n,p},\,\left\langle \boldsymbol{q}_{m,\centerdot,p},\,\boldsymbol{r}_{n,\centerdot,m},\,\boldsymbol{s}_{p,\centerdot,n}\right\rangle -\delta_{m,n,p}\:\right\langle _{1\le m,n,p\le l}
\end{equation}
where $1\le m,n,p\le l$. which corresponds to a subset of the polynomial
ring over the indicated set of variables. The following theorem is
equivalent to theorem 1.\\
\textbf{Theorem: }\emph{(for $3$-tensors)} If $\boldsymbol{A}$ is
a non zero hermitian and $\left\Vert \boldsymbol{A}\right\Vert _{\ell_{3}}^{3}\ne1$
then the spectral system of equations expressed as
\begin{equation}
\begin{cases}
\begin{array}{c}
\left\langle \boldsymbol{D}_{m,p}\cdot\boldsymbol{q}_{m,\centerdot,p},\:\boldsymbol{E}_{n,m}\cdot\boldsymbol{r}_{n,\centerdot,m},\:\boldsymbol{F}_{p,n}\cdot\boldsymbol{s}_{p,\centerdot,n}\right\rangle =a_{m,n,p}\\
\left\langle \boldsymbol{q}_{m,\centerdot,p},\,\boldsymbol{r}_{n,\centerdot,m},\,\boldsymbol{s}_{p,\centerdot,n}\right\rangle =\delta_{m,n,p}
\end{array}\end{cases}
\end{equation}
admits a solution.

\subsubsection*{\emph{Proof: }}

Similarly to the $2$-tensor case, we exhibit a polynomial $p$ which
does not belong to the Ideal $\mathcal{I}$ defined bellow. 
\begin{equation}
\mathcal{I}=\left\rangle \:\left\langle \boldsymbol{D}_{m,p}\cdot\boldsymbol{q}_{m,\centerdot,p},\:\boldsymbol{E}_{n,m}\cdot\boldsymbol{r}_{n,\centerdot,m},\:\boldsymbol{F}_{p,n}\cdot\boldsymbol{s}_{p,\centerdot,n}\right\rangle -a_{m,n,p},\,\left\langle \boldsymbol{q}_{m,\centerdot,p},\,\boldsymbol{r}_{n,\centerdot,m},\,\boldsymbol{s}_{p,\centerdot,n}\right\rangle -\delta_{m,n,p}\:\right\langle _{1\le m\le n\le p\le l}.
\end{equation}
Such a polynomial $p$ is expressed by 
\[
p=\left(\sum_{1\le i,j,k\le l}\left|\left\langle \boldsymbol{D}_{m,p}\cdot\boldsymbol{q}_{m,\centerdot,p},\:\boldsymbol{E}_{n,m}\cdot\boldsymbol{r}_{n,\centerdot,m},\:\boldsymbol{F}_{p,n}\cdot\boldsymbol{s}_{p,\centerdot,n}\right\rangle \right|^{3}\right)^{3}-\left\Vert \boldsymbol{A}\right\Vert _{\ell_{3}}^{3}
\]
\begin{equation}
p\notin\mathcal{I}
\end{equation}
since 
\begin{equation}
p\in\mathcal{I}\Rightarrow\left\Vert \boldsymbol{A}\right\Vert _{\ell_{3}}^{3^{2}}=\left\Vert \boldsymbol{A}\right\Vert _{\ell_{3}}^{3}
\end{equation}
which contradicts our assumption that $\left\Vert \boldsymbol{A}\right\Vert _{\ell_{3}}^{3}\ne1$,
this completes the proof. $\square$\\
Hermiticity also ensure that the solution to the spectral decomposition
is not the trivial all zero solution since for non zero $3$-tensor
$\boldsymbol{A}$ 
\begin{equation}
\left\Vert \boldsymbol{A}\right\Vert _{\ell_{3}}^{3}=\left(\sum_{1\le k\le l}\left\{ \circ\left(\boldsymbol{A},\boldsymbol{A},\boldsymbol{A}\right)\right\} _{k,k,k}+\sum_{1\le i<j<k\le l}a_{i,j,k}\cdot a_{k,i,j}^{\mathfrak{c}_{3}^{2}}\cdot a_{j,k,i}^{\mathfrak{c}_{3}^{1}}\right)>0
\end{equation}

\section{Properties following from the spectral decomposition}

Similarly to the formulation for the spectral theorem for matrices,
we can also discuss the notion of eigen-objects for tensors. In order
to point out the analogy let us consider the matrix decomposition
equations in Eq \ref{eq:matrix_decomp1} and Eq \ref{eq:matrix_decomp2},
one is therefore led to consider the matrices $\widetilde{\boldsymbol{Q}}\equiv\left(\tilde{q}_{m,n}=\sqrt{\lambda_{m}}\, q_{n,m}\right)$
as the scaled matrix of eigenvectors. According to our proposed decomposition,
the corresponding equations for $3$-tensors is given by

\begin{equation}
a_{m,n,p}=\sum_{1\le k\le l}\left(\mu_{m,k}\cdot q_{m,k,p}\cdot\mu_{k,p}\right)\left(\nu_{n,k}\cdot r_{n,k,m}\cdot\nu_{k,m}\right)^{\mathfrak{c}_{3}^{2}}\left(\xi_{p,k}\cdot s_{p,k,n}\cdot\xi_{k,n}\right)^{\mathfrak{c}_{3}^{1}},
\end{equation}
recall that the tensor $\widetilde{\boldsymbol{Q}}:=\left(\tilde{q}_{m,k,p}=\omega_{m,k}\cdot\omega_{k,p}\cdot q_{m,k,p}\right)$
collects as slices what we refer to as the \emph{scaled eigen-matrices}.
The analogy with eigenvectors is based on the following outerproduct
expansion.

\begin{equation}
\boldsymbol{A}=\sum_{1\le k\le l}\left(\left(\mu_{k}\,\boldsymbol{q}_{k,\centerdot}\right)^{\mathfrak{c}_{2}^{1}}\otimes\;\left(\nu_{k}\,\boldsymbol{r}_{k,\centerdot}\right)\right).
\end{equation}
The equation emphasizes the fact that a hermitian matrices can be
viewed as a sum of exterior products of scaled eigenvectors and the
scaling factor associated to the rank one matrix resulting from the
outerproduct corresponds to the eigenvalue. Similarly, a symmetric
$3$-tensor may also be viewed as a sum outer products of slices or
matrices and therefore we refer to the corresponding slices as scaled
\emph{eigen-matrices}. The outerproduct sum follows from the identity

\begin{equation}
\boldsymbol{A}=\circ\left(\widetilde{\boldsymbol{Q}},\widetilde{\boldsymbol{R}}^{\dagger^{2}},\widetilde{\boldsymbol{S}}^{\dagger}\right)
\end{equation}
expressed as :

\begin{equation}
\boldsymbol{A}=\sum_{k=1}^{l}\otimes\left(\tilde{\mathbf{q}}_{\centerdot,k,\centerdot},\tilde{\mathbf{r}}_{\centerdot,\centerdot,k},\tilde{\mathbf{s}}_{k,\centerdot,\centerdot}\right),\label{eq:Outer product expansion}
\end{equation}
 which can be equivalently written as
\[
a_{m,n,p}=
\]

\begin{equation}
\sum_{1\le k\le l}\left(\left(\mu_{m,k}\cdot\mu_{k,p}\right)\left(\nu_{n,k}\cdot\nu_{k,m}\right)^{\mathfrak{c}_{3}^{2}}\left(\xi_{p,k}\cdot\xi_{k,n}\right)^{\mathfrak{c}_{3}^{1}}\right)f_{m,n,p}(k)
\end{equation}
 where $f_{m,n,p}(k)$ denote the $k$-th component expressed 

\begin{equation}
f_{m,n,p}(k)\::=q_{m,k,p}\left(r_{n,k,m}\right)^{\mathfrak{c}_{3}^{2}}\left(s_{p,k,n}\right)^{\mathfrak{c}_{3}^{1}}.
\end{equation}
\emph{We may summarize by simply saying that: as one had eigenvalues
and eigenvectors for matrices one has eigenvectors and eigen-matrices
for $3$-tensors.}

\section{Computational Framework}

We shall first provide an algorithmic description of the characteristic
polynomial of matrix without assuming the definition of the determinant
of matrices and furthermore show how the description allows us to
define characteristic polynomials for tensors. We recall for a matrix
that the characteristic system of equations is determined by the algebraic
system of equations 
\begin{equation}
\boldsymbol{A}=\boldsymbol{Q}^{T}\cdot\boldsymbol{D}\cdot\boldsymbol{R}\Leftrightarrow\begin{cases}
\begin{array}{c}
\left\langle \boldsymbol{D}^{\frac{1}{2}}\cdot\boldsymbol{q}_{m},\boldsymbol{D}^{\frac{1}{2}}\cdot\boldsymbol{r}_{n}\right\rangle =a_{m,n}\\
\left\langle \boldsymbol{q}_{m},\boldsymbol{r}_{n}\right\rangle =\delta_{m,n}
\end{array} & 1\le m\le n\le l\end{cases}
\end{equation}
as discussed above induces the following polynomial ideal

\begin{equation}
\mathcal{I}=\left\rangle \;\left\langle \boldsymbol{D}^{\frac{1}{2}}\cdot\boldsymbol{q}_{m},\boldsymbol{D}^{\frac{1}{2}}\cdot\boldsymbol{r}_{n}\right\rangle -a_{m,n},\:\left\langle \boldsymbol{q}_{m},\boldsymbol{q}_{n}\right\rangle -\delta_{m,n}\:\right\langle _{1\le m\le n\le l}\subseteq\mathbb{C}\left[\left\{ \lambda_{k},\:\boldsymbol{q}_{k},\boldsymbol{r}_{k}\right\} _{1\le k\le l}\right].
\end{equation}
Let $\mathcal{G}$ be the reduced Gr$\ddot{o}$bner basis of $\mathcal{I}$
using the ordering on the monomials induced by the following lexicographic
ordering of the variables. 
\begin{equation}
\boldsymbol{Q}>\boldsymbol{R}>\lambda_{1}>\cdots>\lambda_{l}
\end{equation}
In the case of matrices it has been established that there is a polynomial
relationship between the eigenvalues; more specifically the eigenvalues
are roots to the algebraic equation

\begin{equation}
p(\lambda)=\det\left(\boldsymbol{A}-\lambda\cdot\boldsymbol{I}\right)
\end{equation}
By the elimination theorem \cite{Cox-Little-O'shea} we may computationaly
derive the characteristic polynomials as follows 

\begin{equation}
\mathcal{I}\cap\mathbb{C}\left[\lambda_{l}\right]=\det\left(\boldsymbol{A}-\lambda_{l}\boldsymbol{I}\right)
\end{equation}
It therefore follows from this observation that the reduced Gr$\ddot{o}$bner
basis of $\mathcal{I}$ determines the characteristic polynomial of
$\boldsymbol{A}$.

\subsection*{Definition}

Let $\mathcal{G}$ denote the reduced Gr$\ddot{o}$bner basis of the
ideal $\mathcal{I}$ using the the lexicographic order on the monimials
induced by the following lexicographic order of the variables.
\[
\boldsymbol{Q}>\boldsymbol{R}>\boldsymbol{S}>\boldsymbol{D}>\boldsymbol{E}>\boldsymbol{F}
\]
where 
\[
\mathcal{I}=\left\rangle \:\left\langle \boldsymbol{D}_{m,p}\cdot\boldsymbol{q}_{m,\centerdot,p},\:\boldsymbol{E}_{n,m}\cdot\boldsymbol{r}_{n,\centerdot,m},\:\boldsymbol{F}_{n,m}\cdot\boldsymbol{s}_{p,\centerdot,n}\right\rangle -a_{m,n,p},\,\left\langle \boldsymbol{q}_{m,\centerdot,p},\,\boldsymbol{r}_{n,\centerdot,m},\,\boldsymbol{s}_{p,\centerdot,n}\right\rangle -\delta_{m,n,p}\:\right\langle _{1\le m\le n\le p\le l}
\]
The reduced characteristic set of polynomials $\mathcal{C}$ associated
with the hermitian $3$-tensor $\boldsymbol{A}$ is a subset of the
reduced Groebner basis $\mathcal{G}$ such that 

\begin{equation}
\mathcal{C}\::=\mathcal{G}\cap\mathbb{C}\left[\boldsymbol{D},\boldsymbol{E},\boldsymbol{F}\right]
\end{equation}
where $\mathbb{C}\left[\boldsymbol{D},\boldsymbol{E},\boldsymbol{F}\right]$
denotes the polynomial ring in the entries of the sacaling tensor
with complex coefficients. The reduced should here be thougth of as
generalization of the characteristic polynomial associated with matrices.

\section{The General Framework}

\subsection{$n$-tensor Algebra}

An $\left(m_{1}\times m_{2}\times\cdots\times m_{n-1}\times m_{n}\right)$
$n$-tensor $\boldsymbol{A}$ is a set of elements of a field indexed
by the set resulting from the Cartesian product

\[
\left\{ 1,2,\cdots,(m_{1}-1),m_{1}\right\} \times\left\{ 1,2,\cdots,(m_{2}-1),m_{2}\right\} \times\cdots\times\left\{ 1,2,\cdots,(m_{n}-1),m_{n}\right\} 
\]
The dimensions of $\boldsymbol{A}$ is specified by $\left(m_{1}\times m_{2}\times\cdots\times m_{n-1}\times m_{n}\right)$
where $\forall$ $1\le k\le n$ , $m_{k}\in\mathbb{N}^{\star}$ specifies
the dimensions of the tensor. We may also introduce a dimension operator
defined by

\begin{equation}
d(\boldsymbol{A},k)=\begin{cases}
\begin{array}{cc}
m_{k} & if\;1\le k\le n\\
0 & else
\end{array}\end{cases}
\end{equation}
 Finally, we shall simply use the notation convention $\boldsymbol{A}=\left(a_{i_{1},i_{2},\cdots,i_{n}}\right)$
for describing $\boldsymbol{A}$ once the dimensions have been specified.\\
In what follows we will discuss general tensor products for $n$-tensors
where $n$ is a positive integer greater or equal to 2. Let us start
by recalling the definition of matrix multiplication

\begin{equation}
b_{i_{1},i_{2}}=\sum_{j}a_{i_{1},j}^{(1)}\cdot a_{j,i_{2}}^{(2)},
\end{equation}
the preceding matrix product generalizes to the proposed $3$-tensor
product as follows

\begin{equation}
b_{i_{1},i_{2},i_{3}}=\sum_{j}a_{i_{1},j,i_{3}}^{(1)}\cdot a_{i_{1},i_{2},j}^{(2)}\cdot a_{j,i_{2},i_{3}}^{(3)}.
\end{equation}
By closely inspecting the expression of the product we note that if
\textbf{$\boldsymbol{A}^{(1)}$} is a $\left(m\times k\times1\right)$
tensor, and $\boldsymbol{A}^{(3)}$ is a $\left(k\times n\times1\right)$
tensor then the resulting tensor $\boldsymbol{B}$ expressed by

\begin{equation}
b_{i_{1},i_{2},1}=\sum_{j}a_{i_{1},j,1}^{(1)}\cdot a_{i_{1},i_{2},j}^{(2)}\cdot a_{j,i_{2},1}^{(3)}\:\forall\left(i_{1},i_{2}\right)\: s.t.\;\left(\begin{array}{c}
1\le i_{1}\le m\\
1\le i_{2}\le n
\end{array}\right)\label{eq:Tensor Action}
\end{equation}
 will be of dimensions $\left(m\times n\times1\right)$. The product
above expresses the action of $3$-tensor $\boldsymbol{A}^{(2)}$
of dimension $\left(m\times n\times k\right)$ on the pair of matrices
arising from $\boldsymbol{A}^{(1)}$ and $\boldsymbol{A}^{(3)}$.
Furthermore for $\boldsymbol{A}^{(2)}$ having entries such that

\begin{equation}
\boldsymbol{A}^{(2)}\equiv\left(a_{i_{1},i_{2},j}^{(2)}=1\right)\;\forall\left(i_{1},i_{2},j\right)\: s.t.\;\left(\begin{array}{c}
1\le i_{1}\le m\\
1\le i_{2}\le n\\
1\le j\le k
\end{array}\right),
\end{equation}
the result of the action of $\boldsymbol{A}^{(2)}$ on the pair of
matrices arising from the tensors$\boldsymbol{A}^{(1)}$ and $\boldsymbol{A}^{(3)}$
simply corresponds to a matrix multiplication. For $4$-tensor the
product operator is expressed as :

\begin{equation}
b_{i_{1},i_{2},i_{3},i_{4}}=\sum_{j}a_{i_{1},j,i_{3},i_{4}}^{(1)}\cdot a_{i_{1},i_{2},j,i_{4}}^{(2)}\cdot a_{i_{1},i_{2},i_{3},j}^{(3)}\cdot a_{j,i_{2},i_{3},i_{4}}^{(4)}.
\end{equation}
Similarly the tensor $\boldsymbol{A}^{(3)}$ can be chosen to be all-one
tensor which reduces the product above to the product operation for
$3$-tensors. This nested relationship will also apply to higher order
tensors.\\
We may now write the expression for the product of $n$-tensor. Let
$\left\{ \boldsymbol{A}^{(t)}=\left(a_{i_{1},i_{2},\cdots,i_{n}}^{(t)}\right)\right\} _{1\le t\le n}$
denotes a set of $n$-tensors. The product operator has therefore
$n$ operands and is noted:

\begin{equation}
\boldsymbol{B}=\bigcirc_{t=1}^{n}\left(\boldsymbol{A}^{(t)}\right)
\end{equation}
defined by

\begin{equation}
b_{i_{1},i_{2},\cdots,i_{n}}=\sum_{k}\left(a_{i_{1},k,i_{2},\cdots,i_{n}}^{(1)}\times\cdots\times a_{i_{1},i_{2},\cdots,i_{t},\, k\,,i_{t+2},\cdots,\, i_{n}}^{(t)}\times\cdots\times a_{k\,,i_{2},\cdots,i_{n}}^{(n)}\right)
\end{equation}

\begin{equation}
b_{i_{1},i_{2},\cdots,i_{n}}=\sum_{k}\left(\left(\prod_{t=1}^{n-1}a_{i_{1},i_{2},\cdots,i_{t},\, k\,,i_{t+2},\cdots,\, i_{n}}^{(t)}\right)a_{k\,,i_{2},\cdots,i_{n}}^{(n)}\right)
\end{equation}
 It follows from the definition that the dimensions of the tensors
in the set $\left\{ \boldsymbol{A}^{(t)}=\left(a_{i_{1},i_{2},\cdots,i_{n}}^{(t)}\right)\right\} _{1\le t\le n}$must
be chosen so that :

\begin{equation}
d(\boldsymbol{A}^{(1)},2)=d(\boldsymbol{A}^{(2)},3)=\cdots=d(\boldsymbol{A}^{(n-1)},n)=d(\boldsymbol{A}^{(n)},1).
\end{equation}
 which describes the constraints on the dimension relating all the
$n$ tensors in the product. The constraints accross the $(n-1)$
other dimensions for each tensor are described by the following relation.

\begin{equation}
d(\boldsymbol{A}^{(i)},k)=d(\boldsymbol{A}^{(j)},k)\;\forall k\notin\left\{ \left(j+1\right),\:\left(i+1\right)\right\} 
\end{equation}
 The tensor $\boldsymbol{B}$ resulting from the product is a $n$-tensor
of dimensions .

\begin{equation}
\left(d\left(\boldsymbol{A}^{(1)},1\right)\times d\left(\boldsymbol{A}^{(2)},2\right)\times\cdots\times d\left(\boldsymbol{A}^{(n-1)},(n-1)\right)\times d\left(\boldsymbol{A}^{(n)},n\right)\right)
\end{equation}
 \emph{Note that the product of tensors of lower order all arise as
special cases of the general product formula describe above.}

\emph{Tensor Action:}\\
The action of $n^{th}$ order tensor $\boldsymbol{A}=\left(a_{i_{1},i_{2},\cdots,i_{n}}\right)$
on $\left(n-1\right)$-tuple of order $\left(n-1\right)$ tensors
$\left\{ \boldsymbol{B}^{(t)}=\left(b_{1,i_{2},\cdots,i_{n}}^{(t)}\right)\right\} _{1\le t\le\left(n-1\right)}$is
defined as

\begin{equation}
b_{1,i_{2},\cdots,i_{n}}=\sum_{k}\left(\left(\prod_{t=1}^{n-1}b_{1,i_{2},\cdots,i_{t},\, k\,,i_{t+2},\cdots,\, i_{n}}^{(t)}\right)a_{k\,,i_{2},\cdots,i_{n}}^{(n)}\right).
\end{equation}
 The equation above generalizes the notion of matrices action on a
vector.

\emph{Tensor Outerproduct: }The outer-product of $n$-tuple $\left(n-1\right)$-tensors
is denoted by :

\begin{equation}
\boldsymbol{B}=\bigotimes_{t=1}^{n}\left(\boldsymbol{A}^{(t)}\right)
\end{equation}
 and defined such that :

\begin{equation}
b_{i_{1},i_{2},\cdots,i_{n}}=\left(\left(\prod_{t=1}^{n-1}a_{i_{1},i_{2},\cdots,i_{t},\,1\,,i_{t+2},\cdots,\, i_{n}}^{(t)}\right)a_{1\,,i_{2},\cdots,i_{n}}^{(n)}\right).
\end{equation}
The Kronecker $n$-tensor is defined as

\begin{equation}
\boldsymbol{\Delta}=\left(\delta_{i_{1},i_{2},\cdots,i_{(n-1)},i_{n}}=\left(\prod_{t=1}^{n-1}\delta_{i_{t},i_{(t+1)}}\right)\delta_{i_{n},i_{1}}\right)\equiv\sum_{k}\left(\vec{\boldsymbol{e}}_{k}^{\otimes n}\right)
\end{equation}
\emph{Order $n$ tensor transpose/adjoint: }\\
Given a tensor $\boldsymbol{A}=\left(a_{j_{1},j_{2},\cdots,j_{n}}\right)$
the transpose $\boldsymbol{A}^{T}$ is defined such that

\begin{equation}
\boldsymbol{A}^{T}=\left(a_{j_{2},j_{3},\cdots,j_{n},j_{1}}\right).
\end{equation}
 For a complex valued tensor where the entries are expressed in their
polar form as follows :

\begin{equation}
\boldsymbol{A}=\left(a_{j_{1},j_{2},\cdots,j_{n}}=r_{j_{1},j_{2},\cdots,j_{n}}\cdot exp\left\{ i\cdot\theta_{j_{1},j_{2},\cdots,j_{n}}\right\} \right),
\end{equation}
 the generalized adjoint is given by

\begin{equation}
\boldsymbol{A}^{\dagger}=\left(r_{j_{2},j_{3},\cdots,j_{n},j_{1}}\cdot exp\left\{ i\cdot\exp\left\{ i\cdot\frac{2\pi}{n}\right\} \cdot\theta_{j_{2},j_{3},\cdots,j_{n},j_{1}}\right\} \right),
\end{equation}
\begin{equation}
\boldsymbol{A}^{\dagger^{k}}=\left(r_{\sigma_{k}\left(j_{1}\right),\sigma_{k}\left(j_{2}\right),\cdots,\left(j_{n}\right)}\cdot exp\left\{ i\cdot\exp\left\{ i\cdot\frac{2\pi k}{n}\right\} \cdot\theta_{j_{2},j_{3},\cdots,j_{n},j_{1}}\right\} \right),
\end{equation}
where $\sigma_{k}$ denotes the composition of $k$ cyclic permutation
of the indices from which it follows that 

\begin{equation}
\boldsymbol{A}^{\dagger^{n}}=\boldsymbol{A}.
\end{equation}

\subsection{The Spectrum of $n$-tensors.}

In order to formulate the spectral theorem for $\boldsymbol{A}\in\mathbb{C}^{l^{n}}$
we will briefly discussed notion of orthogonal and scaling $n$-tensors,
which can be expressed as

\begin{equation}
\boldsymbol{\Delta}=\bigcirc_{t=1}^{n}\left(\boldsymbol{Q}^{\dagger^{(n+1-t)}}\right)
\end{equation}
 that is

\begin{equation}
\delta_{i_{1},i_{2},\cdots,i_{n}}=\sum_{k}\left(\left(\prod_{t=1}^{n-1}q_{i_{1},i_{2},\cdots,i_{t},\, k\,,i_{t+2}\cdots,\, i_{n}}^{\dagger^{(n+1-t)}}\right)q_{k,\, i_{2},\cdots,i_{n}}^{\dagger}\right),
\end{equation}
 Where $T$ denotes the transpose operation, which still corresponds
to a cyclic permutation of the indices.

We first provide the formula for the scaling tensor whose product
with $\boldsymbol{A}$ leaves the tensor unchanged.

\begin{equation}
a_{i_{1},i_{2},\cdots,i_{n}}=\left(\bigcirc\left(\boldsymbol{A},\boldsymbol{D}^{(1)},\boldsymbol{D}^{(2)},\boldsymbol{D}^{(3)},\cdots,\boldsymbol{D}^{(n-1)}\right)\right)_{i_{1},i_{2},\cdots,i_{n}}
\end{equation}

\begin{equation}
\Rightarrow a_{i_{1},i_{2},\cdots,i_{n}}=\sum_{k}\left(a_{i_{1},k,i_{2},\cdots,i_{n}}\times d_{i_{1},i_{2},k,\cdots,i_{n}}^{(1)}\times\cdots\times d_{i_{1},i_{2},\cdots,i_{t},\, k\,,i_{t+2},\cdots,\, i_{n}}^{(t)}\times\cdots\times d_{k\,,i_{2},\cdots,i_{n}}^{(n-1)}\right)
\end{equation}

\begin{equation}
\Rightarrow\begin{cases}
\begin{array}{c}
\forall t<n-2\quad\boldsymbol{D}^{(t)}\equiv\left(d_{i_{1},i_{2},\cdots,i_{n}}^{(t)}=\delta_{i_{2},i_{2+t}}\right)\\
\boldsymbol{D}^{(n-1)}\equiv\left(d_{i_{1},i_{2},\cdots,i_{n}}^{(n-1)}=\delta_{i_{1},i_{2}}\right)
\end{array}\end{cases}
\end{equation}
The above family of tensors play the role of identity operator and
are related to one another by transposition of the indices. The more
general expression for the scaling tensors is therefore given by

\begin{equation}
\begin{cases}
\begin{array}{c}
\forall t<n-2\quad\boldsymbol{S}^{(t)}\equiv\left(s_{i_{1},i_{2},\cdots,i_{n}}^{(t)}=\delta_{i_{2},i_{2+t}}\cdot\omega_{i_{t},i_{2+t}}\right)\\
\boldsymbol{S}^{(n-1)}\equiv\left(s_{i_{1},i_{2},\cdots,i_{n}}^{(n-1)}=\delta_{i_{1},i_{2}}\cdot\omega_{i_{1},i_{n-1}}\right)
\end{array}\end{cases}
\end{equation}
where $\boldsymbol{W}=\left(w_{m,n}\right)$ is a symmetric matrix.
The expression for the scaled orthogonal tensor is therefore expressed
by

\begin{equation}
\left(\bigcirc\left(\boldsymbol{Q},\boldsymbol{S}^{(1)},\boldsymbol{S}^{(2)},\boldsymbol{S}^{(3)},\cdots,\boldsymbol{S}^{(n-1)}\right)\right)_{i_{1},i_{2},\cdots,i_{n}}=q_{i_{1},i_{2},\cdots,i_{n}}\left(\prod_{k\ne2}\omega_{i_{2},i_{k}}\right)
\end{equation}
We therefore obtain that the scaled tensor which will be of the form
:

\begin{equation}
\widetilde{\boldsymbol{Q}}=\bigcirc\left(\boldsymbol{Q},\boldsymbol{S}^{(1)},\boldsymbol{S}^{(2)},\boldsymbol{S}^{(3)},\cdots,\boldsymbol{S}^{(n-1)}\right)
\end{equation}

\noindent \textbf{Theorem 2:}\emph{ (Spectral Theorem for $n$-Tensors):}
For any non zero hermitian tensor $\boldsymbol{A}\in\mathbb{C}^{l^{n}}$
such that $\left\Vert \boldsymbol{A}\right\Vert _{\ell_{n}}^{n}\ne1$,
there exist a factorization in the form

\begin{equation}
\begin{cases}
\begin{array}{c}
\boldsymbol{A}=\bigcirc_{t=1}^{n}\left(\widetilde{\boldsymbol{Q}}_{t}^{\dagger^{(n+1-t)}}\right)\\
\boldsymbol{\Delta}=\bigcirc_{t=1}^{n}\left(\boldsymbol{Q}_{t}^{\dagger^{(n+1-t)}}\right)
\end{array}\end{cases}
\end{equation}
the expression above generalizes Eq \ref{eq:Tensor factorization scaled}

\subsubsection*{Proof of the Spectral Theorem for $n$-tensors}

The spectral decompostion yields the following system of equations 

\begin{equation}
\begin{cases}
\begin{array}{c}
\boldsymbol{A}=\bigcirc_{t=1}^{n}\left(\widetilde{\boldsymbol{Q}}_{t}^{\dagger^{(n+1-t)}}\right)\\
\boldsymbol{\Delta}=\bigcirc_{t=1}^{n}\left(\boldsymbol{Q}_{t}^{\dagger^{(n+1-t)}}\right)
\end{array}\end{cases}
\end{equation}
more insightfully rewritten as 

\begin{equation}
\begin{cases}
\begin{array}{c}
\left\langle \boldsymbol{D}_{i_{1},i_{3},\cdots,\, i_{n}}^{(1)}\cdot\boldsymbol{q}_{i_{1},\centerdot,i_{3},\cdots,\, i_{n}}^{(1)},\cdots,\:\boldsymbol{D}_{i_{1},\cdots,i_{t},i_{t+2}\cdots,\, i_{n}}^{(t)}\cdot\boldsymbol{q}_{i_{1},\cdots,i_{t},\,\centerdot\,,i_{t+2}\cdots,\, i_{n}}^{(t)},\cdots,\boldsymbol{D}_{i_{2},\cdots,i_{n}}^{(n)}\cdot\boldsymbol{q}_{\centerdot,\, i_{2},\cdots,i_{n}}^{(n)}\right\rangle =a_{i_{1},i_{2},\cdots,i_{n}}\\
\left\langle \boldsymbol{q}_{i_{1},\centerdot,i_{3},\cdots,\, i_{n}}^{(1)},\cdots,\boldsymbol{q}_{i_{1},\cdots,i_{t},\,\centerdot\,,i_{t+2}\cdots,\, i_{n}}^{(t)},\cdots,\boldsymbol{q}_{\centerdot,\, i_{2},\cdots,i_{n}}^{(t)}\right\rangle =\delta_{i_{1},i_{2},\cdots,i_{n}}
\end{array}\end{cases}
\end{equation}
where $\boldsymbol{D}_{i_{1},\cdots,i_{t},i_{t+2}\cdots,\, i_{n}}^{(t)}$
is a diagonal matrix whose entries are specified by 
\begin{equation}
\boldsymbol{D}_{i_{1},\cdots,i_{t},i_{t+2}\cdots,\, i_{n}}^{(t)}=\left(d_{m,n}^{i_{1},\cdots,i_{t},i_{t+2}\cdots,\, i_{n}}(t)=\delta_{m,n}\cdot\omega_{m,n}\right)
\end{equation}
We had already pointed out earlier in the proof for the spectral theorem
for $3$-tensors that the proof technique would apply to $n$-tensors
with norm $\ne1$, where $n$ is a positive integer greater or equal
to $2$. Similarly we consider the polynomial expression
\[
p=\left(\sum_{1\le i_{1},\cdots,i_{n}\le l}\right.
\]
\[
\left.\left|\left\langle \boldsymbol{D}_{i_{1},i_{3},\cdots,\, i_{n}}^{(1)}\cdot\boldsymbol{q}_{i_{1},\centerdot,i_{3},\cdots,\, i_{n}}^{(1)},\cdots,\:\boldsymbol{D}_{i_{1},\cdots,i_{t},i_{t+2}\cdots,\, i_{n}}^{(t)}\cdot\boldsymbol{q}_{i_{1},\cdots,i_{t},\,\centerdot\,,i_{t+2}\cdots,\, i_{n}}^{(t)},\cdots,\boldsymbol{D}_{i_{2},\cdots,i_{n}}^{(t)}\cdot\boldsymbol{q}_{\centerdot,\, i_{2},\cdots,i_{n}}^{(n)}\right\rangle \right|^{n}\right)^{n}
\]
 
\begin{equation}
-\left\Vert \boldsymbol{A}\right\Vert _{\ell_{n}}^{n}
\end{equation}
and observe that 
\begin{equation}
p\notin\mathcal{I}
\end{equation}
where $\mathcal{I}$ defines the ideal iduced by the spectral system
of equation since 
\begin{equation}
p\in\mathcal{I}\Rightarrow\left\Vert \boldsymbol{A}\right\Vert _{\ell_{n}}^{n^{2}}=\left\Vert \boldsymbol{A}\right\Vert _{\ell_{n}}^{n}
\end{equation}
which contradicts our assumption that $\left\Vert \boldsymbol{A}\right\Vert _{\ell_{n}}^{n}\ne1$,
Hence we conclude that 
\begin{equation}
\left\Vert \boldsymbol{A}\right\Vert _{\ell_{n}}^{n}\ne1\Rightarrow p\notin\mathcal{I}
\end{equation}
this completes the proof. $\square$ \\
The $l$ {}``slices'' of the scaled tensor $\widetilde{\boldsymbol{Q}}_{t}$
constitutes what we call the \emph{scaled eigen-tensors} of $\boldsymbol{A}$
which are $\left(n-1\right)$-tensors.

\subsection{Spectral Hierarchy}

We recursively define the spectral hierarchy for a tensor $\boldsymbol{A}\in\mathbb{C}^{l^{n}}$
. The base case for the recursion is the case of matrices. The spectrum
of an $\left(l\times l\right)$ matrix is characterized by a set of
$l$ scaled eigen-vectors. The existence of the spectral hierarchy
relies on the observation that the spectrum of an order $n$-tensor
$\boldsymbol{A}\in\mathbb{C}^{l^{n}}$ is determined by a collection
of $l$-tuple $\left(n-1\right)$-tensors not necessarily distinct.
Each one of these $l$-tuples corresponding to a scaled orthogonal
eigen-tensor. By recursively computing the spectrum of the resulting
scaled orthogonal $\left(n-1\right)$-tensors, one determines a tree
structure which completely characterizes the spectral hierarchy associated
with the $n$-tensor $\boldsymbol{A}$. The leaves of the tree will
be made of scaled eigenvectors when the spectral decomposition exists
for all the resulting lower order tensors.

It therefore follows that the tensor $\boldsymbol{A}$ can be expressed
as a nested sequence of sums of outer products. We illustrate the
general principle with $3$-tensors. Let $\boldsymbol{A}$ denotes
a third order tensor which admits a spectral decomposition in the
form described by Eq \ref{eq:Outer product expansion}. We recall
that the spectral decomposition for $3$-tensors is expressed by

\begin{equation}
\boldsymbol{A}=\circ\left(\tilde{\boldsymbol{Q}},\,\tilde{\boldsymbol{R}}^{\dagger^{2}},\,\tilde{\boldsymbol{S}}^{\dagger}\right)
\end{equation}

\begin{equation}
\boldsymbol{A}=\sum_{k=1}^{l}\otimes\left(\left(\mu_{m,k}\cdot\mu_{k,p}\cdot q_{m,k,p}\right)_{m,p},\left(\nu_{n,k}\cdot\nu_{k,m}\cdot r_{n,k,m}\right)_{n,m},\left(\xi_{p,k}\cdot\xi_{k,n}\cdot r_{p,k,n}\right)_{p,n}\right)
\end{equation}
 by computing the spectrum of the scaled eigen-matrices we have :

\begin{equation}
\forall\:1\le j_{1}\le l\quad\tilde{\boldsymbol{Q}}(k)=\left(\mu_{m,k}\cdot\mu_{k,p}\cdot q_{m,k,p}\right)_{m,p}=\sum_{1\le j_{1}\le l}\left(\sqrt{\gamma_{j_{1}}(k)}\cdot\vec{\boldsymbol{u}}_{j_{1}}(k)\right)\otimes\left(\sqrt{\gamma_{j_{1}}(k)}\cdot\vec{\boldsymbol{v}}_{j_{1}}(k)\right)
\end{equation}

\begin{equation}
\forall\:1\le j_{2}\le l\quad\tilde{\boldsymbol{R}}(k)=\left(\nu_{n,k}\cdot\nu_{k,m}\cdot r_{n,k,m}\right)_{n,m}=\sum_{1\le j_{2}\le l}\left(\sqrt{\lambda_{j_{2}}(k)}\cdot\vec{\boldsymbol{w}}_{j_{2}}(k)\right)\otimes\left(\sqrt{\lambda_{j_{2}}(k)}\cdot\vec{\boldsymbol{x}}_{j_{2}}(k)\right)
\end{equation}

\begin{equation}
\forall\:1\le j_{3}\le l\quad\tilde{\boldsymbol{S}}(k)=\left(\xi_{p,k}\cdot\xi_{k,n}\cdot r_{p,k,n}\right)_{n,m}=\sum_{1\le j_{3}\le l}\left(\sqrt{\beta_{j_{3}}(k)}\cdot\vec{\boldsymbol{y}}_{j_{3}}(k)\right)\otimes\left(\sqrt{\beta_{j_{3}}(k)}\cdot\vec{\boldsymbol{z}}_{j_{3}}(k)\right)
\end{equation}
 where $\forall\;1\le k\le l$ , $\gamma_{j_{1}}(k),\left\{ \vec{\boldsymbol{u}}_{j_{1}}(k),\vec{\boldsymbol{v}}_{j_{1}}(k)\right\} $,$\lambda_{j_{2}}(k)$,$\left\{ \vec{\boldsymbol{w}}_{j_{2}}(k),\vec{\boldsymbol{x}}_{j_{2}}(k)\right\} $
and $\beta_{j_{3}}(k)$,$\left\{ \vec{\boldsymbol{y}}_{j_{3}}(k),\vec{\boldsymbol{z}}_{j_{3}}(k)\right\} $
denote the eigenvalues and corresponding eigenvectors respectively
for the matrices $\tilde{\boldsymbol{S}}(k)$ ,$\tilde{\boldsymbol{Q}}(k)$,
$\tilde{\boldsymbol{R}}(k)$. It therefore follows that $\boldsymbol{A}$
can be expressed by the following nested sum of outer product expressions

\[
\boldsymbol{A}=
\]

\[
\sum_{k=1}^{l}\otimes\left(\left[\sum_{1\le j_{2}\le l}\left(\sqrt{\gamma_{j_{1}}(k)}\cdot\vec{\boldsymbol{u}}_{j_{1}}(k)\right)\otimes\left(\sqrt{\gamma_{j_{1}}(k)}\cdot\vec{\boldsymbol{v}}_{j_{1}}(k)\right)\right],\right.
\]

\begin{equation}
\left.\left[\sum_{1\le j_{2}\le l}\left(\sqrt{\lambda_{j_{2}}(k)}\cdot\vec{\boldsymbol{w}}_{j_{2}}(k)\right)\otimes\left(\sqrt{\lambda_{j_{2}}(k)}\cdot\vec{\boldsymbol{x}}_{j_{2}}(k)\right)\right],\left[\sum_{1\le j_{3}\le l}\left(\sqrt{\beta_{j_{3}}(k)}\cdot\vec{\boldsymbol{y}}_{j_{3}}(k)\right)\otimes\left(\sqrt{\beta_{j_{3}}(k)}\cdot\vec{\boldsymbol{z}}_{j_{3}}(k)\right)\right]\right)
\end{equation}

\section{Relation to previously proposed decompositions}

We shall present in this section a brief overview of the relationship
between our framework and earlier proposed tensor decompositions

\subsection{Tucker Decomposition.}

Let us show in this section how the Tucker decomposition in fact uses
matrix algebra more specifically orthogonality of matrices to express
the singular value decomposition for $3$-tensors. We used for this
section the notation and convention we introduced through this work.
The Tucker factorization scheme finds for an arbitrary $3$-tensor
$\boldsymbol{D}$ the following decomposition

\begin{equation}
\boldsymbol{D}=\boldsymbol{T}\times_{1}\boldsymbol{Q}^{(1)}\times_{2}\boldsymbol{S}^{(2)}\times_{3}\boldsymbol{U}^{(3)},
\end{equation}
 where $\boldsymbol{T}$ denotes a $3$-tensor and $\boldsymbol{Q}^{(1)},\boldsymbol{S}^{(2)},\boldsymbol{U}^{(3)}$
denote matrices. The product expression used for the decomposition
written above corresponds to our proposed definition for triplet dot
product with non trivial background as described in Eq \ref{eq:Matrices to Tensor}.
Using our notation we can express the decomposition of $\boldsymbol{D}$
as follows:

\begin{equation}
d_{m,n,p}=\langle a_{m,i,1},b_{1,n,j},c_{k,1,p}\rangle_{\boldsymbol{T}}=\sum_{i}\sum_{j}\sum_{k}a_{m,i,1}\cdot b_{1,n,j}\cdot c_{k,1,p}\cdot t_{i,j,k}
\end{equation}
 Our starting point is the following invariance relation, which arises
from the matrix products with the identity matrix.

\begin{equation}
d_{m,n,p}=\sum_{i}\sum_{j}\sum_{k}\gamma_{m,i,1}\cdot\gamma_{1,n,j}\cdot\gamma_{k,1,p}\cdot d_{i,j,k}\:,\label{eq:Invariance Tucker}
\end{equation}
 where $\gamma_{m,i,1}=\delta_{m,i}$ , $\gamma_{1,n,j}=\delta_{n,j}$
and $\gamma_{k,1,p}=\delta_{k,p}$ which correspond to transposes
of the identity matrix. For any orthogonal matrices $\boldsymbol{Q}$,
$\boldsymbol{S}$ and $\boldsymbol{U}$ we know that

\begin{equation}
\begin{cases}
\begin{array}{c}
\gamma_{m,i,1}=\sum_{y}q_{m,y,1}\cdot q_{i,y,1}\\
\gamma_{1,n,j}=\sum_{r}s_{1,n,r}\cdot s_{1,j,r}\\
\gamma_{k,1,p}=\sum_{v}u_{k,1,v}\cdot u_{p,1,v}
\end{array}\end{cases}
\end{equation}
Incidentally the expression in Eq \ref{eq:Invariance Tucker} can
be written as :

\begin{equation}
\sum_{i}\sum_{j}\sum_{k}\left(\sum_{y}q_{m,y,1}\cdot q_{i,y,1}\right)\cdot\left(\sum_{r}s_{1,n,r}\cdot s_{1,j,r}\right)\cdot\left(\sum_{v}u_{k,1,v}\cdot u_{p,1,v}\right)\cdot d_{i,j,k}
\end{equation}
by interchanging the order of the sums we get :

\begin{equation}
\sum_{y}\sum_{r}\sum_{v}\left(\sum_{i}q_{m,y,1}q_{i,y,1}\right)\cdot\left(\sum_{j}s_{1,n,r}s_{1,j,r}\right)\cdot\left(\sum_{k}u_{k,1,v}\cdot u_{p,1,v}\right)d_{i,j,k}
\end{equation}
we now separate out the products in the expressions to yield the general
form of the Tucker decomposition.

\begin{equation}
\Rightarrow\sum_{y}\sum_{r}\sum_{v}q_{m,y,1}\cdot s_{1,n,r}\cdot u_{p,1,v}\left(\sum_{i}\sum_{j}\sum_{k}q_{i,y,1}\cdot s_{1,j,r}\cdot u_{k,1,v}\cdot d_{i,j,k}\right)
\end{equation}

\begin{equation}
T\equiv\left(t_{y,r,v}=\sum_{i}\sum_{j}\sum_{k}q_{i,y,1}\cdot s_{1,j,r}\cdot u_{k,1,v}\cdot d_{i,j,k}\right)
\end{equation}
The preceding emphasizes that the Tucker decomposition reuses matrix
orthogonality and does not provide a generalization of the notion
of orthogonality for $n$-tensors. Finally to determine the orthogonal
matrices $\boldsymbol{Q}$, $\boldsymbol{S}$ and $\boldsymbol{U}$
to use we specify the following constraints

\begin{equation}
\sum_{l}\sum_{g}t_{l,g,\alpha}\cdot t_{l,g,\beta}=\delta_{\alpha,\beta}\cdot\left(\sum_{l,g}\left(t_{l,g,\alpha}\right)^{2}\right)
\end{equation}

\begin{equation}
\sum_{l}\sum_{g}t_{l,\alpha,g}\cdot t_{l,\beta,g}=\delta_{\alpha,\beta}\cdot\left(\sum_{l,g}\left(t_{l,\alpha,g}\right)^{2}\right)
\end{equation}

\begin{equation}
\sum_{l}\sum_{g}t_{\alpha,l,g}\cdot t_{\beta,l,g}=\delta_{\alpha,\beta}\cdot\left(\sum_{l,g}\left(t_{\alpha,l,g}\right)^{2}\right)
\end{equation}
which is referred to as the total orthogonality condition.

\subsection{Tensor Rank 1 decomposition.}

The Rank 1 decomposition of tensor \cite{Qi_L3,Grigoriev98,Hastad90,Raz2010,Hillar2009,Grigoriev78,Grigoriev81}
corresponds to solving the following optimization problem. Given an
$r$-tensor $\boldsymbol{A}=\left(a_{i_{1},\cdots,i_{r}}\right)$
we seek to find:

\begin{equation}
\min_{\left(\boldsymbol{x}_{k}^{(t)}\right)_{1\le t\le r}\in\left(\bigotimes_{1\le t\le r}V_{t}\right)}||\boldsymbol{A}-\sum_{1\le k\le l}\left(\lambda_{k}\right)^{r}\bigotimes_{1\le t\le r}\vec{\boldsymbol{x}}_{k}^{(t)}||\label{eq:Tensor rank}
\end{equation}
Since Johan H$\ddot{a}$stad in \cite{J.Hastsad}established the intractability
of the tensor rank problem for $3$-tensors we briefly discuss the
relationship to our framework. It follows from the definition of the
outer product of matrices to form a $3$-tensor that

\begin{equation}
\otimes\left(\boldsymbol{M}_{1}\equiv\left(m_{s,1,t}\right)_{s,t},\boldsymbol{N}_{1}\equiv\left(n_{s,t,1}\right)_{s,t},\boldsymbol{P}_{1}\equiv\left(p_{1,s,t}\right)_{s,t}\right)\equiv\boldsymbol{D}\equiv\left(d_{i,j,k}=m_{i,1,k}\cdot n_{i,j,1}\cdot p_{1,j,k}\right).
\end{equation}
 We point out that for the very special matrices essentially made
up of the same vector as depicted bellow :

\begin{equation}
m_{i,1,k}=u_{i,1,1}\;\forall\:1\le k\le l
\end{equation}

\begin{equation}
n_{i,j,1}=v_{1,j,1}\;\forall\:1\le i\le l
\end{equation}

\begin{equation}
p_{1,j,k}=w_{1,1,k}\;\forall\:1\le j\le l
\end{equation}
 the outer product of the matrices

\begin{equation}
\otimes\left(\boldsymbol{M}_{1}\equiv\left(m_{s,1,t}\right)_{s,t},\boldsymbol{N}_{1}\equiv\left(n_{s,t,1}\right)_{s,t},\boldsymbol{P}_{1}\equiv\left(p_{1,s,t}\right)_{s,t}\right)=\vec{\boldsymbol{u}}\otimes\vec{\boldsymbol{v}}\otimes\vec{\boldsymbol{w}}.
\end{equation}
 This allows us to formulate the tensor rank problem in Eq \ref{eq:Tensor rank}
in terms of the outer product operator for slices as follows

\begin{equation}
\min\left\Vert \left(\sum_{1\le k\le l}\otimes\left(\boldsymbol{M}_{k}\equiv\left(\lambda_{k}\cdot m_{s,k,t}\right)_{s,t},\boldsymbol{N}_{k}\equiv\left(\lambda_{k}\cdot n_{s,t,k}\right)_{s,t},\boldsymbol{P}_{k}\equiv\left(\lambda_{k}\cdot p_{k,s,t}\right)_{s,t}\right)\right)-\boldsymbol{A}\right\Vert _{\ell_{3}}
\end{equation}

\begin{equation}
\Leftrightarrow\min\left\Vert \circ\left(\boldsymbol{M},\boldsymbol{N},\boldsymbol{P}\right)-\boldsymbol{A}\right\Vert _{\ell_{3}},
\end{equation}
 where $\boldsymbol{M},\boldsymbol{N},\boldsymbol{P}$ are $3$-tensors
arising from the collection of matrices associated with the collection
of vectors. The preceding naturally related the tensor rank problem
to our proposed tensor product. Furthermore the generalized framework
allows us to formulate the tensor rank problem for $n$-tensor where
$n$ is a positive integer greater or equal to 2 as follows

\begin{equation}
\min\left\Vert \left(\bigcirc_{t=1}^{n}\left(\boldsymbol{M}^{(t)}\right)\right)-\boldsymbol{A}\right\Vert _{\ell_{n}}\label{eq:Tensor Rank Problem}
\end{equation}
 One may point out that the spectral decomposition associated with
a Hermitian tensor comes quite close to the sought after decomposition
at the cost of the trading of the requirement that the matrices should
be rank one to the fact the matrices should arise from scaled eigen-tensors.

\section{Conclusion}

In this paper we introduced a generalization of the spectral theory
for $n$-tensors where $n$ is a positive integer greater or equal
to $2$. We propose a mathematical framework for $3$-tensors algebra
based on a ternary product operator, which generalizes to $n$-tensors.
This algebra allows us to generalize notions and operators we are
familiar with from Linear algebra including dot product, tensor adjoints,
tensor hermicity, diagonal tensor, permutation tensors and characteristic
polynomials. We proved the spectral theorem for tensors having Forbenius
norm different from $1$. Finally we discussed the spectral hierarchy
which confirms the intractability of determining the orthogonal vector
components whose exterior product result in a given $n$-tensor.

Starting from the recently proposed product formula in Eq \ref{eq:ternaryproduct}
for order $3$-tensors proposed by P. Bhattacharya in \cite{Bhat}
we were able to formulate a general algebra for finite order tensors.
The order $3$-tensor product formula suggests a definition for outer
product of matrices as discussed in Eq \ref{eq:Outer-Product}, it
also suggests how to express the action of a tensor on lower order
tensors. Most importantly with Eq \ref{eq:Background-tensors} we
propose a natural generalization for the dot product operator and
a generalization for the Riemann metric tensor ideas. Furthermore
the tensor algebra that we discuss sketches possible approaches to
investigate generalizations of inner product space theory.

One important characteristic of the product operator for tensor of
order strictly greater than 2 is the fact that the product is not
associative. Incidentally by analogy to matrix theory where the lost
of commutativity for matrix product results into a commutator theory
and lie Algeras which plays an important role in quantum mechanics,
the lost of associativity as expressed in Eq \ref{eq:Non_Associative}
could potentially give rise to an associator theory or generalizations
of lie algebras. Furthermore the transpose operator described in Eq
\ref{eq:Transpose} emphasizes the importance of the roots of unity
in generalizing herminian and unitary tensors. The $3$-tensor permutation
tensors provided a suprising representation for the permutation group
$S_{n}$ which provide a glimpse at a tensor approach to a representation
theory as well as a tensor approach to Markov tensor models.

At the heart of our work lies the concept of orthogonal tensors. We
emphasize the fact the orthogonal tensors discussed here are generalizations
of orthogonal matrices and are significantly different from orthogonal
matrices. One significant difference lie in the two distinct interpretation
of the orthogonality property for tensor. The first interpretation
expressed by Eq \ref{eq:Orthonormaility} is analogous to orthonormal
for a set of vectors. The second interpretation relates to the invariance
of the Kronecker delta tensor under conjugation as expressed in Eq\ref{eq:Kronecker Invariance}.
Furthermore we have through this work provided a natural generalization
for the familiar characteristic polynomial using the important tool
set of Grobner Basis.

Spectral analysis plays an important role in the theory and investigations
of Graphs. Graph spectra have proved to be a relatively useful graph
invariant for determining Isomorphism class of graphs. It seem of
interest to note that the symmetries of a graph described by it's
corresponding automorphism group can also be viewed as depicting a
3-uniform hypergraph which can in turn be investigated by through
it spectral properties. Determining the relationship between spectral
properties of a graph and the spectral properties of it corresponding
automorphism seems worthy of attention in the context of determining
isomorphism classes of graphs. The general framework which address
the algebra for arbitrarily finite order tensor allowed us to derive
the spectral hierarchy. The spectral hierarchy induces a bottom up
construction for finite order tensor from vectors. This explicit construction
may in fact prove useful in the context investigations on tensor rank
problems which also validate as illustrated in Eq \ref{eq:Tensor Rank Problem}
our product operator.

\subsection*{Acknowledgment:}

We are grateful to Emilie Hogan, Professor Doron Zeilberger and Professor
Henry Cohn for helpful discussion regarding properties of Ideals.
The first author was partially supported by the National Science Foundation
grant NSF-DGE-0549115.

\pagebreak{} 
\end{document}